\documentclass[11pt]{article}
\usepackage{amsmath}
\usepackage{amsthm}
\usepackage{epsfig}
\usepackage{leftidx}
\usepackage{graphicx,subfigure}
\usepackage{amsfonts}
\usepackage[ruled,boxed,commentsnumbered,linesnumbered]{algorithm2e}
\usepackage{graphicx, color, epstopdf}
\usepackage{cite}
\usepackage{flafter}
\usepackage{amssymb,amsmath}
\usepackage{multicol}
\usepackage{booktabs}
\usepackage{multirow}
\usepackage{enumerate}
\usepackage{enumitem}
\usepackage{algpseudocode,algorithmicx}
\usepackage{xcolor}
\usepackage{tikz}
\usetikzlibrary{arrows,shapes,chains}

\newtheorem{theorem}{Theorem}[section]
\newtheorem{lemma}{Lemma}[section]
\newtheorem{example}{Example}
\newtheorem{remark}{Remark}
\newtheorem{definition}{Definition}

\newcommand{\rmq}{\mathrm{q}}
\newcommand{\rmk}{\mathrm{k}}
\newcommand{\two}{\mathbf{2}}
\newcommand{\three}{\mathbf{3}}
\newcommand{\four}{\mathbf{4}}
\newcommand{\piff}{\partial_\tau}
\newcommand{\defeq}{:=}
\newcommand{\zd}{\,\mathrm{d}}
\newcommand{\abs}[1]{\left|#1\right|}
\newcommand{\absb}[1]{\big|#1\big|}

\newcommand{\bra}[1]{\left(#1\right)}
\newcommand{\brab}[1]{\big(#1\big)}
\newcommand{\braB}[1]{\Big(#1\Big)}
\newcommand{\brat}[1]{(#1)}
\newcommand{\kbra}[1]{\left[#1\right]}

\newcommand{\kbraB}[1]{\Big[#1\Big]}

\title{Mesh-robust stability and convergence of variable-step deferred correction methods based on the  BDF2 formula}
\author{\quad Jiahe Yue \thanks{School of Mathematics, Nanjing University of Aeronautics and Astronautics,
		Nanjing 211106, China. Email: yuejiahe@nuaa.edu.cn.}
	\quad
	Hong-lin Liao\thanks{ORCID 0000-0003-0777-6832. School  of Mathematics, 
		Nanjing University of Aeronautics and Astronautics,
		Nanjing 211106, China; Key Laboratory of Mathematical Modeling
		and High Performance Computing of Air Vehicles (NUAA), MIIT, Nanjing 211106, China. 
		Emails: liaohl@nuaa.edu.cn and liaohl@csrc.ac.cn.
		This author's work is supported by NSF of China
		under grant number 12071216.}
	\quad
	Nan Liu\thanks{School  of Mathematics, Nanjing University of Aeronautics and Astronautics,
				Nanjing 211106, China. Email: liunan@nuaa.edu.cn.}				
	}
\date{\today}

\begin{document}
	
	\maketitle
	
	\begin{abstract}
		We provide a new theoretical framework for the variable-step deferred correction (DC) methods based on the well-known BDF2 formula. By using the discrete orthogonal convolution kernels, some high-order BDF2-DC methods are proven to be stable on arbitrary time grids according to the recent definition of stability (SINUM, 60: 2253-2272). It significantly relaxes the existing step-ratio restrictions for the BDF2-DC methods (BIT, 62: 1789-1822). The associated sharp error estimates are established by taking the numerical effects of the starting approximations into account, and they suggest that the BDF2-DC methods have no aftereffect, that is, the lower-order starting scheme for the BDF2 scheme will not cause a loss in the accuracy of the high-order BDF2-DC methods. Extensive tests on the graded and random time meshes are presented to support the new theory.
		\\[1ex]
		\emph{Keywords}: variable-step deferred correction methods, high order methods,
		discrete orthogonal convolution kernels, stability and convergence
		\\[1ex]
		\emph{AMS subject classifications}: 65M06, 65M12
	\end{abstract}
	
	\section{Introduction}
	In many applications like acoustics, electromagnetics and fluid dynamics, 
	stable high-order numerical methods \cite{AkrivisFeischlKovacsLubich2021,AlbiDimarcoPareschi:2020,BouchritiPierreAlaa2020,Couture-PeckDelfour:2020,CrouzeixLisbona:1984,DimarcoPareschi:2017,DuttGreengardRokhkin2000,GustafssonKress2001,BourgaultGaron:2022,Butcher1987,Gear1971,GongZhaoWang2020,HairerNrsettWanner1993,HairerWanner1996,Lambert1991,LiaoKang2023,LiaoTangZhou2023BDF3parabolic} are always required. Roughly speaking, the existing methods can be divided into two types: intrinsically high-order schemes \cite{BouchritiPierreAlaa2020,Butcher1987,Couture-PeckDelfour:2020,CrouzeixLisbona:1984,DimarcoPareschi:2017,Gear1971,KangLiaoWang:2023cnsns,HairerNrsettWanner1993,HairerWanner1996,LiaoKang2023,LiaoTangZhou2023BDF3parabolic,Lambert1991,GongZhaoWang2020,AkrivisFeischlKovacsLubich2021} and methods based on Richardson extrapolation or deferred correction \cite{DuttGreengardRokhkin2000,GustafssonKress2001,BourgaultGaron:2022} to accelerate the low-order schemes. Usually,  the high-order schemes tend to have relatively poor stability properties, especially when they are applied to simulate multi-scale physical phenomena, such as the chaotic motion of submerged devices \cite{Couture-PeckDelfour:2020}, the hyperbolic systems with multiscale relaxation \cite{AlbiDimarcoPareschi:2020}
	and the stiff kinetic equations \cite{DimarcoPareschi:2017}. Actually, the difficulty of choosing the optimal time-step size to capture the time-sensitive variations of physical phenomena limits the predictive efficiency of these time integration methods. In such applications, the variable-step time-stepping methods, such as the variable-step BDF formulas \cite{Gear1971,CrouzeixLisbona:1984,KangLiaoWang:2023cnsns,LiaoKang2023,LiLiao:2022,LiaoJiZhang:2020pfc,LiaoTangZhou2023,LiaoZhang:2021,LiaoTangZhou:2020bdf2,LiaoTangZhou2023BDF3parabolic}, are preferred since they allow larger (smaller) time-step sizes when the physical solution is slowly varying (rapidly changing), cf. \cite{BouchritiPierreAlaa2020,HairerNrsettWanner1993,LiaoJiZhang:2020pfc,LiaoTangZhou:2020bdf2}.
	
	Here we consider the variable-step deferred correction (DC) methods \cite{BourgaultGaron:2022} based on the second-order backward differentiation formula (BDF2) for the first-order differential problems
	\begin{align}\label{con: ODE equation}
		\frac{\zd v}{\zd t}=f(t,v)\quad\text{for $0<t\le T$}\quad \text{with $v(0)=v_0$},
	\end{align}
    where the nonlinear function $f$ is Lipschitz-continuous.
	Without loss of generality, put time levels $0=t_0<t_1<t_2<\cdots<t_N=T$ with the variable time-step
	$\tau_k\defeq t_k-t_{k-1}$ for $1\le k\le N$, and the maximum step size
	$\tau\defeq\max_{1\le k\le N}\tau_k$. Define the adjacent time-step ratios
	$r_{k}\defeq \tau_{k}/\tau_{k-1}$ for $k\ge2$ with $r_1:=0$, and denote $r_{\max}:=\max_{k\ge2}r_k$.
	We consider a time-discrete solution, $v^n\approx v(t_n)$, of the following
	 BDF2 scheme
	\begin{align}\label{eq: BDF2 scheme}
		D_2 v^n=f(t_n,v^n),
	\end{align}
	where the BDF2 operator $D_2$ is defined by
	\begin{align*}
		D_2v^n:=\sum_{j=1}^nd^{(n)}_{n-j}\piff  v^j
		\quad\text{for $n\ge2$}
	\end{align*}
	with  $\piff v^n:=(v^n-v^{n-1})/\tau_n$ and the discrete convolution kernels 
	$d^{(n)}_{n-j}$ are defined as
	\begin{align*}
		&d^{(n)}_0:=\frac{1+2r_n}{1+r_n},\quad d^{(n)}_1:=-\frac{r_n}{1+r_n}\quad\text{and}\quad
		d^{(n)}_j:=0\quad\text{for $n\ge j+1\ge 2$}.
	\end{align*}
	By substituting the exact solution $v(t)$ into the BDF2 scheme \eqref{eq: BDF2 scheme}, one has
	\begin{align}\label{con: ODE DC equation}
		D_2 v(t_n)+R^n(t_n)=f(t_n,v(t_n))\quad\text{for $n\ge 2$},
	\end{align}
	where the truncation error 
	$$R^n(t_n)=\sum_{j=3}^\infty(-1)^{j+1}\frac{v^{(j)}(t_n)}{j!}\bra{-(1+r_n)\tau_n^{j-1}
		+r_n(\tau_n+\tau_{n-1})^{j-1}}.$$
	
	To obtain higher accuracy, the key point of the BDF2-based DC (BDF2-DC) methods \cite{BourgaultGaron:2022} is to approximate the truncation error $R^n(t_n)$ by a $\rmq$-order($\rmq=3,4$) deferred correction $\mathcal{D}_{2,\rmq}f(t_n,v^n)$ that can be written as 
	\begin{align}\label{def: DC2q operator}
		\mathcal{D}_{2,\rmq} f(t_n,v^n):=\sum_{j=3}^\rmq (-1)^{j+1}\frac{v^{(j),n}}{j!}\brab{-(1+r_n)\tau_n^{j-1}+r_n(\tau_n+\tau_{n-1})^{j-1}}
	\end{align}
	for $n\geq \rmq-1$, where $v^{(j),n}$ is the numerical approximation of $v^{(j)}(t_n)$. These DC methods with Taylor's remainder originated from the construction of the difference corrected BDF framework in \cite{Soderlind1989}. Subsequently, Gustafsson \cite{GustafssonKress2001} applied them to the trapezoidal rule (called the Crank–Nicolson scheme when applied to partial differential problems) and obtained the $\rmq$-order methods based on the deferred correction principle.
	
	Specifically, using the Newton difference quotient $f[\cdot,\cdot,\ldots,\cdot,\cdot]$, cf.  \cite{BourgaultGaron:2022}, one has the following quadratic Newton interpolation at points $t_n,\,t_{n-1}$ and $t_{n-2}$ ($n\geq 2$),
	\begin{align*}
		N_2(t)&=f(t_n,v(t_n))+f[t_n,t_{n-1}](t-t_n)+f[t_n,t_{n-1},t_{n-2}](t-t_n)(t-t_{n-1}).
	\end{align*} 
	Then, taking $\rmq=3$ in \eqref{def: DC2q operator} and using the quadratic Newton interpolation $N_2(t_n)$ to approximate $f(t_n,v(t_n))$, it is not difficult to get
	\begin{align}\label{de: D23}
		\mathcal{D}_{2,3} f(t_n,v^n) := \frac{v^{(3),n}}{3!}\kbra{-(1+r_n)\tau_n^{2}+r_n(\tau_n+\tau_{n-1})^{2}},
	\end{align} 
	where $v^{(3),n}:=\frac{d^2N_2(t_n)}{dt^2}=2f[t_n,t_{n-1},t_{n-2}]$.
	Similarly, considering the following cubic Newton interpolation 
	at points $t_n,\,t_{n-1},\,t_{n-2}$ and $t_{n-3}$ ($n\geq 3$),
	it is easy to find 
	\begin{align*}
		N_3(t)&=f(t_n,v^n)+f[t_n,t_{n-1}](t-t_n)+f[t_n,t_{n-1},t_{n-2}](t-t_n)(t-t_{n-1})\\
		&\quad+f[t_n,t_{n-1},t_{n-2},t_{n-3}](t-t_n)(t-t_{n-1})(t-t_{n-2}).
	\end{align*} 
	Then, by taking $\rmq=4$ in \eqref{def: DC2q operator} and employing the cubic Newton interpolation $N_3(t_n)$ to approximate $f(t_n,v(t_n))$, we have
	\begin{align}\label{de: D24}
		\mathcal{D}_{2,4} f(t_n,v^n) :&= \frac{v^{(3),n}}{3!}\kbra{r_n(\tau_n+\tau_{n-1})^{2}-(1+r_n)\tau_n^{2}}
		-\frac{v^{(4),n}}{4!}\kbra{r_n(\tau_n+\tau_{n-1})^{3}-(1+r_n)\tau_n^{3}}\nonumber\\
		&=\mathcal{D}_{2,3}\,f(t_n,v^n)
		+\frac{\tau_n}{12}(\tau_n+\tau_{n-1})(2\tau_n+\tau_{n-1})f[t_n,t_{n-1},t_{n-2},t_{n-3}]\nonumber\\
		&\triangleq \mathcal{D}_{2,3}\,f(t_n,v^n)+\mathcal{D}_{3,4}\,f(t_n,v^n),
	\end{align} 
	where
	\begin{align*}
		&v^{(4),n}:=\,\frac{d^3N_3(t_n)}{dt^3}=6f[t_n,t_{n-1},t_{n-2},t_{n-3}],\\
		&v^{(3),n}:=\,\frac{d^2N_3(t_n)}{dt^2}=2f[t_n,t_{n-1},t_{n-2}]
		+2f[t_n,t_{n-1},t_{n-2},t_{n-3}](2\tau_n+\tau_{n-1}),\\
		&\mathcal{D}_{3,4}\,f(t_n,v^n):=
		\frac{\tau_n}{12}(\tau_n+\tau_{n-1})(2\tau_n+\tau_{n-1})f[t_n,t_{n-1},t_{n-2},t_{n-3}].
	\end{align*}

	The BDF2-DC methods are step-by-step solution procedures \cite{BourgaultGaron:2022}. 
	Firstly, one uses the BDF2 scheme to calculate the second-order solution $v_\two ^n$. Then we use the third-order correction $\mathcal{D}_{2,3} f(t_n,v_\two ^n)$ to construct the BDF2-DC3 scheme, which generates the third-order solution $v_\three^n$. To achieve fourth-order accurate solution $v_\four^{n}$, the BDF2-DC3-DC4 scheme is constructed by adding the fourth-order correction $\mathcal{D}_{2,4} f(t_n,v_\three^n)$ in the BDF2 code. This series of algorithms leads to the following numerical schemes
	\begin{align}
		\textbf{BDF2:}& \quad
		D_2v_\two ^{n} = f(t_n,v_\two ^n)\quad \text{for $n\geq 2$},\label{eq: DC2 scheme}\\
		\textbf{BDF2-DC3:}&\quad
		D_2 v_\three^{n} + \mathcal{D}_{2,3} f(t_n,v_\two ^n)  
		= f(t_n,v_\three^n)\quad\text{for $n\geq 2$},\label{eq: BDF2-DC3 scheme}\\
		\textbf{BDF2-DC3-DC4:}&\quad 
		D_2 v_\four^{n} + \mathcal{D}_{2,4} f(t_n,v_\three^n)   
		= f(t_n,v_\four^n)\quad\text{for $n\geq 3$}. \label{eq: BDF2-DC3-DC4 scheme}	
	\end{align}

	\begin{figure}[htb!]
		\centering
		\scriptsize
		\tikzstyle{format}=[rectangle,draw,thin,fill=white]
		\begin{tikzpicture}
			\node (BDF2){\textbf{BDF2:}};
			\node[format,right of=BDF2,node distance=20mm] (v20){$v_\two ^0$};
			\node[format,right of=v20,node distance=25mm] (v21){$v_\two ^1$};
			\node[format,right of=v21,node distance=25mm] (v22){$v_\two ^2$};
			\node[format,right of=v22,node distance=25mm] (v23){$v_\two ^3$};
			\node[right of=v23,node distance=25mm] (v24){$\cdots$};
			
			\node[below of=BDF2,node distance=20mm] (BDF2-DC3){\textbf{BDF2-DC3:}};
			\node[format,right of=BDF2-DC3,node distance=20mm] (v30){$v_\three^0$};
			\node[format,right of=v30,node distance=25mm] (v31){$v_\three^1$};
			\node[format,right of=v31,node distance=25mm] (v32){$v_\three^2$};
			\node[format,right of=v32,node distance=25mm] (v33){$v_\three^3$};
			\node[right of=v33,node distance=25mm] (v34){$\cdots$};	
			
			\draw[->] (v20)--node[above]{BDF1/RK}(v21);
			\draw[->](v21)--node[above]{(1.7)}(v22);
			\draw[->](v22)--node[above]{(1.7)}(v23);
			\draw[->](v23)--node[above]{(1.7)}(v24);
			\draw[->] (v22)--node[left]{$\mathcal{D}_{2,3}$}(v32);
			\draw[->] (v23)--node[left]{$\mathcal{D}_{2,3}$}(v33);
			
			\draw[->] (v30)--node[above]{BDF1/RK}(v31);
			\draw[->](v31)--node[above]{(1.8)}(v32);
			\draw[->](v32)--node[above]{(1.8)}(v33);
			\draw[->](v33)--node[above]{(1.8)}(v34);
		\end{tikzpicture}
		\caption{The starting procedure (BDF1 or RK) 
			of BDF2-DC3 scheme \eqref{eq: BDF2-DC3 scheme}.}\label{figure: DC3 initial process}
	\end{figure}
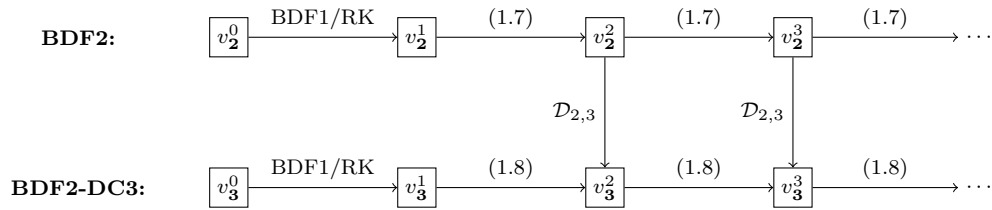

\begin{figure}[htb!]
		\centering
		\scriptsize
		\tikzstyle{format}=[rectangle,draw,thin,fill=white]
		\begin{tikzpicture}		
			\node[below of=BDF2-DC3,node distance=20mm] (BDF2-DC33){\textbf{BDF2-DC3:}};
			\node[format,right of=BDF2-DC33,node distance=20mm] (v330){$v_\three^0$};
			\node[format,right of=v330,node distance=25mm] (v331){$v_\three^1$};
			\node[format,right of=v331,node distance=25mm] (v332){$v_\three^2$};
			\node[format,right of=v332,node distance=25mm] (v333){$v_\three^3$};
			\node[right of=v333,node distance=25mm] (v334){$\cdots$};
			
			\node[below of=BDF2-DC33,node distance=20mm] (BDF2-DC3-DC4){\textbf{BDF2-DC3-DC4:}};
			\node[format,right of=BDF2-DC3-DC4,node distance=20mm] (v40){$v_\four^0$};
			\node[format,right of=v40,node distance=25mm] (v41){$v_\four^1$};
			\node[format,right of=v41,node distance=25mm] (v42){$v_\four^2$};
			\node[format,right of=v42,node distance=25mm] (v43){$v_\four^3$};
			\node[right of=v43,node distance=25mm] (v44){$\cdots$};	
			
			\draw[->] (v330)--node[above]{BDF1/RK}(v331);
			\draw[->](v331)--node[above]{(1.8)}(v332);
			\draw[->](v332)--node[above]{(1.8)}(v333);
			\draw[->](v333)--node[above]{(1.8)}(v334);
			\draw[->] (v333)--node[left]{$\mathcal{D}_{2,4}$}(v43);
			
			\draw[->] (v40)--node[above]{BDF1/RK}(v41);
			\draw[->](v41)--node[above]{BDF1/RK}(v42);
			\draw[->](v42)--node[above]{(1.9)}(v43);
			\draw[->](v43)--node[above]{(1.9)}(v44);
		\end{tikzpicture}
		\caption{The starting procedure (BDF1 or RK) of 
			BDF2-DC3-DC4 scheme \eqref{eq: BDF2-DC3-DC4 scheme}.}
		\label{figure: DC4 initial process}
	\end{figure}
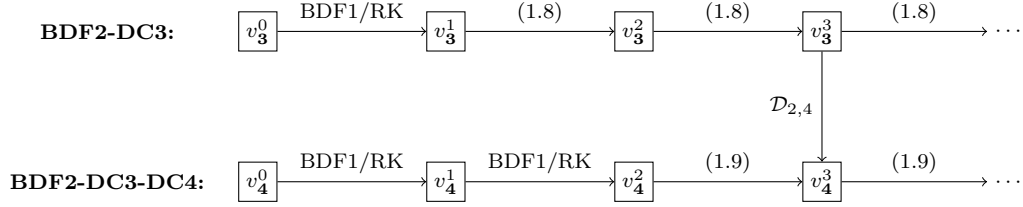

	Since $D_2$ and $\mathcal{D}_{2,3}$ are three-level operators according to \eqref{de: D23}, we need the starting values $v_\two ^1$ and $v_\three^1$ at the first level $t_1$ to startup the BDF2 scheme \eqref{eq: DC2 scheme} and the  BDF2-DC3 scheme \eqref{eq: BDF2-DC3 scheme}, see Figure \ref{figure: DC3 initial process}. Meanwhile, since the operator $\mathcal{D}_{2,4}$ involves the solutions at four levels according to the definition \eqref{de: D24}, two starting values $v_\four^1$ and $v_\four^2$ are necessary to startup the  BDF2-DC3-DC4 scheme \eqref{eq: BDF2-DC3-DC4 scheme},  see Figure \ref{figure: DC4 initial process}. Obviously, different resolutions of these starting solutions $v_\two ^1$, $v_\three^1$, $v_\four^1$ and $v_\four^2$ would arrive at different accuracies of numerical solutions. We always assume that these starting values are available and accurate in establishing our theory while, in the numerical tests, some of one-step approaches including the BDF1 scheme and Runge-Kutta (RK) methods with different orders are examined for potential applications. 

	For the subsequent stability analysis, assume that the perturbed solution $\bar{v}_\two ^{n},\bar{v}_\three^{n}$ and $\bar{v}_\four^{n}$ solve the following perturbed equations, respectively,
	\begin{align*}
		\textbf{BDF2:}& \quad
		D_2\bar{v}_\two ^{n} = f(t_n,\bar{v}_\two ^n)+\varepsilon_\two ^n\quad \text{for $n\geq 2$},\\
		\textbf{BDF2-DC3:}&\quad
		D_2 \bar{v}_\three^{n} + \mathcal{D}_{2,3} f(t_n,\bar{v}_\two ^n)= f(t_n,\bar{v}_\three^n)+\varepsilon_\three^n\quad \text{for $n\geq 2$},\\
		\textbf{BDF2-DC3-DC4:}&\quad 
		D_2 \bar{v}_\four^{n} + \mathcal{D}_{2,4} f(t_n,\bar{v}_\three^n)= f(t_n,\bar{v}_\four^n)+\varepsilon_\four^n\quad \text{for $n\geq 3$},
	\end{align*}
	where the time sequences $\{\varepsilon_\two ^n\},\{\varepsilon_\three^n\}$ and $\{\varepsilon_\four^n\}$ are assumed to be bounded. Let the perturbed errors $\tilde{v}_{\rmq}^n:=\bar{v}_{\rmq}^n-v_{\rmq}^n$ for $\rmq=\mathbf{2},\mathbf{3}$ and $\mathbf{4}$. One has the following perturbed error systems
	\begin{align}
		\textbf{BDF2:}& \quad
		D_2\tilde{v}_\two ^{n} = f(t_n,\bar{v}_\two ^n)-f(t_n,v_\two ^n)+\varepsilon_\two ^n,\label{eq: BDF2 perturbed equation}\\
		\textbf{BDF2-DC3:}&\quad
		D_2 \tilde{v}_\three^{n} = f(t_n,\bar{v}_\three^n)-f(t_n,v_\three^n)+\mathcal{D}_{2,3} f(t_n,v_\two ^n)- \mathcal{D}_{2,3} f(t_n,\bar{v}_\two ^n)+\varepsilon_\three^n,\label{eq: BDF2-DC3 perturbed equation}\\
		\textbf{BDF2-DC3-DC4:}&\quad 
		D_2 \tilde{v}_\four^{n}= f(t_n,\bar{v}_\four^n)-f(t_n,v_\four^n)+\mathcal{D}_{2,4} f(t_n,v_\three^n)-\mathcal{D}_{2,4} f(t_n,\bar{v}_\three^n)+\varepsilon_\four^n.\label{eq: BDF2-DC3-DC4 perturbed equation}
	\end{align}

	To process the analysis of the variable-step multi-step schemes, we recall two different definitions of stability for variable-step $\rmk$-step methods. Compared with the classical zero-stability (Definition \ref{Def:stability1}),  Definition \ref{Def:stability2} includes the discrete time derivatives of starting solutions. 

	\begin{definition}[\textbf{Zero stability}]\label{Def:stability1}\cite[Definition 2.2]{CrouzeixLisbona:1984}
		Assume that the nonlinear function $f(t,v)$ is Lipschitz-continuous with the Lipschitz constant $L_f>0$. A variable-step $\rmk$-step method is stable if there exists a real number $\tau_0$ and a fixed constant $C$ (independent of the step sizes $\tau_n$) such that the perturbed error $\tilde{v}^n:=\bar{v}^n-v^n$ fulfills 
		\begin{align*}
			\max_{\rmk\le i\le N }\abs{\tilde{v}^{i}}\le C\bra{\sum_{i=0}^{\rmk-1}\abs{\tilde{v}^{i}}
				+\max_{\rmk\le i\le N}\abs{\varepsilon^i}}\quad\text{for all $\tau_n\le \tau_0$.}
		\end{align*}
	\end{definition}
	\begin{definition}[\textbf{Stability}]\label{Def:stability2}\cite[Definition 1]{LiLiao:2022}
		Assume that the nonlinear function $f(t,v)$ is Lipschitz-continuous
		with the Lipschitz constant $L_f>0$. A variable-step $\rmk$-step method
		is stable if there exists a real number $\tau_0$ and a fixed constant $C$ (independent of the step sizes $\tau_n$)
		such that	the perturbed error	$\tilde{v}^n:=\bar{v}^n-v^n$  fulfills 
		\begin{align*}
			\max_{\rmk\le i\le N }\abs{\tilde{v}^{i}}\le C\bra{\sum_{i=0}^{\rmk-1}\abs{\tilde{v}^{i}}+\sum_{i=1}^{\rmk-1}\tau\abs{\partial_\tau\tilde{v}^{i}} +\max_{\rmk\le i\le N}\abs{\varepsilon^i}}\quad\text{for all $\tau_n\le \tau_0$.}
		\end{align*}
	\end{definition}
	
	According to Definition \ref{Def:stability1}, Bourgault and Garon \cite{BourgaultGaron:2022} showed that the above BDF2-DC schemes are stable if the adjacent time-step ratios $0<r_k<1+\sqrt{2}$ for $k\ge2$, which is also the classical zero-stability condition for the variable-step BDF2 scheme. Recently, Li and Liao \cite{LiLiao:2022} suggested a new concept of stability for variable-step  multistep schemes, cf. Definition \ref{Def:stability2}, and proved that the BDF2 scheme is  stable on arbitrary time meshes for the problem \eqref{con: ODE equation}. This step-ratio-independent stability was derived with the help of a new analysis tool, namely the discrete orthogonal convolution (DOC) kernels \cite{LiLiao:2022,LiaoJiZhang:2020pfc,LiaoJiWangZhang:2021,LiaoKang2023,LiaoTangZhou2023,LiaoTangZhou2023BDF3parabolic,LiaoTangZhou2023SCM,LiaoZhang:2021}, which was proven to be powerful in the numerical analysis for variable-step discrete time approximations including the BDF schemes and discrete fractional derivatives. In this paper, we will use the DOC kernels to investigate the variable-step BDF2-DC3, BDF2-DC3-DC4 and BDF2-DC4 schemes. Specially, the numerical effects (accuracy) of the starting approximations  $v_\two ^1$, $v_\three^1$, $v_\four^1$ and $v_\four^2$ on these DC solutions are  carefully explored in the analysis and are always reflected in our error estimates.	

	It is to remark that different definitions will lead to different theoretical results. Since we adopt Definition \ref{Def:stability2} in the numerical analysis, the resulting stability and error estimates are significantly different from those in \cite{BourgaultGaron:2022}. The existing step-ratio conditions \cite{BourgaultGaron:2022} for the stability and convergence of 
	the BDF2-DC schemes  \eqref{eq: DC2 scheme}-\eqref{eq: BDF2-DC3-DC4 scheme} are significantly relaxed. Sections 2 and 3 present the stability and error estimates for the  BDF2-DC3 and BDF2-DC3-DC4 schemes, respectively. Section 4 extends our analysis to the variable-step BDF2-DC4 method, which aims to improve two-order of accuracy by one-step deferred correction. Numerical experiments on the graded and random time meshes are included in Section 5 to show that our theoretical results are also practically reliable. 
	
	Throughout the paper, $C$ represents a generic positive constant, not necessarily the same at different occurrences. The appeared constant may be dependent on the given data, the solution and the length of the integration interval but is always independent of the time-step sizes $\tau_n$ and the time-step ratios $r_n$. Specially, a stability or convergence estimate is also called mesh-robust if the associated stability or convergence prefactor (constant) is independent of the step ratios $r_n$. 	
	
    \section{Stability and convergence of BDF2-DC3 scheme}
	\setcounter{equation}{0}
	Assume that the summation $\sum_{k=i}^{j}\cdot$ to be zero and the product $\prod_{k=i}^{j}\cdot$ to be one if the index $i>j$.
	As for the BDF2 kernels $d^{(n)}_{n-j}$ with any fixed indexes $n$, a class of DOC kern $\big\{\vartheta_{n-j}^{(n)}\big\}_{j=2}^n$ are introduced by a recursive procedure \cite{LiaoZhang:2021,LiaoJiWangZhang:2021,LiaoTangZhou2023SCM}
	\begin{align}\label{eq: BDF2-DOC procedure}
		\vartheta_{0}^{(n)}:=\frac{1}{d^{(n)}_{0}}\quad\text{and}\quad
		\vartheta_{n-j}^{(n)}:=-\frac{1}{d^{(j)}_{0}}
		\sum_{i=j+1}^{n}\vartheta_{n-i}^{(n)}d^{(i)}_{i-j}\quad\text{for $2\leq j\le n-1$.}
	\end{align}
	The above DOC kernels $\vartheta_{n-j}^{(n)}$ satisfy the following discrete orthogonality identity
	\cite{LiaoJiWangZhang:2021,LiaoTangZhou2023SCM}
	\begin{align}\label{eq: BDF2-DOC identity}
		\sum_{i=j}^{n}\vartheta_{n-i}^{(n)}d^{(i)}_{i-j}\equiv\delta_{nj},\quad\sum_{i=j}^{n}d_{n-i}^{(n)}\vartheta^{(i)}_{i-j}\equiv\delta_{nj}
		\quad\text{for $2\leq j\le n$,}
	\end{align}
	where $\delta_{nj}$ is the Kronecker delta symbol with $\delta_{nj}=0$ if $j\neq n$. Therefore, multiplying the left of the perturbed error equation \eqref{eq: BDF2 perturbed equation} by the DOC kernels $\vartheta_{n-i}^{(n)}$, summing $i$ from $2$ to $n$ and exchanging the summation order and using \eqref{eq: BDF2-DOC identity}, one has
	\begin{align}\label{eq: BDF2-DOC transform}
		\sum_{i=2}^{n}\vartheta_{n-i}^{(n)}D_{2}\tilde{v}_\two ^i
		=&\,\sum_{i=2}^{n}\vartheta_{n-i}^{(n)}
		d^{(i)}_{i-1}\piff \tilde{v}_\two ^1
		+\sum_{i=2}^{n}\vartheta_{n-i}^{(n)}
		\sum_{j=2}^id^{(i)}_{i-j}\piff \tilde{v}_\two ^j\nonumber\\
		=&\,\piff \tilde{v}_\two ^1
		\sum_{i=2}^n\vartheta_{n-i}^{(n)}d^{(i)}_{i-1}+
		\sum_{j=2}^{n}\piff \tilde{v}_\two ^j
		\sum_{i=j}^n\vartheta_{n-i}^{(n)}d^{(i)}_{i-j}\nonumber\\
		=&\,\mathcal{I}_{1}^n[\tilde{v}_\two ]+\piff   \tilde{v}_\two ^n\qquad\text{for $2\le n\le N$,}
	\end{align}
	where $\mathcal{I}_{1}^n[v]$ represents the starting effect on the numerical solution at the time $t_n$,
	\begin{align}\label{eq: initial effect BDFk-DOC}
		\mathcal{I}_{1}^n[v]:=&\,\piff v^1
		\sum_{i=2}^n\vartheta_{n-i}^{(n)}d^{(i)}_{i-1}=\vartheta_{n-2}^{(n)}d^{(2)}_{1}(\piff v^1)\qquad\text{for $2\le n\le N$.}
	\end{align}
	Multiplying both sides of the perturbed error equation \eqref{eq: BDF2 perturbed equation} by the DOC kernels
	$\vartheta_{n-i}^{(n)}$, and summing $i$ from $2$ to $n$, we apply \eqref{eq: BDF2-DOC transform} to get 
	\begin{align*}
		\piff \tilde{v}_\two ^n
		=-\mathcal{I}_{1}^n[\tilde{v}_\two ]+\sum_{i=2}^n\vartheta_{n-i}^{(n)}
		\kbra{f(t_i,\bar{v}_\two ^i)-f(t_i,v_\two ^i)+\varepsilon_\two ^i}
		\quad\text{for $2\le n\le N$.}
	\end{align*}
   At first, we recall the property of DOC kernels and the stability of BDF2 scheme.
	\begin{lemma}\label{lem: BDF2 orthogonal formula}\cite[Lemma 2.1]{LiLiao:2022}
		The DOC kernels $\vartheta_{n-j}^{(n)}$ in \eqref{eq: BDF2-DOC procedure} have an explicit formula
		\begin{align*}
			\vartheta_{n-j}^{(n)}=\frac{1}{d^{(j)}_{0}}
			\prod_{i=j+1}^n\frac{r_{i}}{1+2r_{i}}>0\quad\text{ for $2\leq j\le n$},
		\end{align*}
		and satisfy 
		\begin{align*}
			\sum_{j=2}^{n}\vartheta_{n-j}^{(n)}=1-\prod_{i=2}^n\frac{r_{i}}{1+2r_{i}}<1\quad\text{for $n\ge2$.}
		\end{align*}
	\end{lemma}
	\begin{lemma}\cite[Theorem 2.1]{LiLiao:2022}\label{lem: BDF2 stability}
		If $\tau\le 1/(4L_f)$, the solution of \eqref{eq: BDF2 perturbed equation}  satisfies
		\begin{align*}
			\abs{\tilde{v}_\two ^{n}}\le 4\exp(4L_ft_{n-1})\braB{\abs{\tilde{v}_\two ^{1}}
				+\tau\abs{\piff \tilde{v}_\two ^{1}}
				+t_{n}\max_{2\le i\le n}\abs{\varepsilon_\two ^i}}\quad\text{for $2\le n\le N$.}
		\end{align*}
		Thus the variable-step BDF2 scheme \eqref{eq: DC2 scheme} is mesh-robustly stable.
	\end{lemma}
	\begin{remark}\label{rem: zero-stable condition}
		According to Definition \ref{Def:stability2}, Lemma \ref{lem: BDF2 stability} states that the variable-step BDF2 scheme is stable for any step-ratio $r_n>0$. But it would not meet with Definition \ref{Def:stability1}.  On the graded time meshes $t_k=T(k/N)^\gamma$ for some $\gamma>1$, the stability estimate relies on the initial data, $$\abs{\tilde{v}_\two ^{1}}+\tau\abs{\piff \tilde{v}_\two ^{1}}\leq \frac{\tau}{\tau_1}\abs{\tilde{v}_\two ^{0}}+\frac{\tau+\tau_1}{\tau_1}\abs{\tilde{v}_\two ^{1}}\leq \frac{2\tau}{\tau_1}\big(\abs{\tilde{v}_\two ^{0}}+\abs{\tilde{v}_\two ^{1}}\big).$$
		Since the stability constant involves a time-step dependent factor $\tau/\tau_1$, one can not call the BDF2 scheme zero-stable according to Definition 1. To be more clear, Section 4 provides some numerical tests, see Tables \ref{table: graded meshes for example1} and \ref{table: graded meshes for example2}, where the solutions are robustly stable and convergent with the desired accuracy; however, the value of the factor $\tau/\tau_1$ increases dramatically and becomes very huge (up to $10^6$) as the step size decreases. In stark contrast, the solution errors decrease gradually as the step size decreases, as expected. In this sense, Definition 2 and the stability result in Lemma \ref{lem: BDF2 stability} (and also Theorems \ref{thm: BDF2-DC3 stability}, \ref{thm: BDF2-DC3-DC4 stability} and \ref{thm: BDF2-DC4 stability} below) seem  practically more reasonable.
	\end{remark}

	\subsection{ Stability of the BDF2-DC3 scheme}
	Now we consider the stability of BDF2-DC3 scheme \eqref{eq: BDF2-DC3 scheme} via its convolutional form. Multiplying both sides of the perturbed error equation \eqref{eq: BDF2-DC3 perturbed equation} by the DOC kernels
	$\vartheta_{n-i}^{(n)}$, and summing $i$ from 2 to $n$, we apply the equality \eqref{eq: BDF2-DOC transform} to get
	\begin{align}\label{eq: BDF2-DC3 DOC perturbed equation}
		\piff \tilde{v}_\three^n
		=-\mathcal{I}_{1}^n[\tilde{v}_\three]+\sum\limits_{i=2}^n\vartheta_{n-i}^{(n)}
		\bra{f(t_i,\bar{v}_\three^i)-f(t_i,v_\three^i)+\mathcal{D}_{2,3}f(t_i,v_\two ^i)-\mathcal{D}_{2,3}f(t_i,\bar{v}_\two ^i)+\varepsilon_\three^i}
	\end{align}
	for $2\le n\le N$. With the definition \eqref{de: D23}, simple calculations arrive at
	\begin{align*}
		\mathcal{D}_{2,3}\,f(t_i,v_\two ^i)
		:&=\frac{1}{3}f[t_i,t_{i-1},t_{i-2}]\big(-(r_i+1)\tau_i^2+\tau_{i-1}\tau_i(r_i+1)^2\big)\\
		&=\frac{\piff  f(t_i,v_\two ^i) - \piff  f(t_{i-1},v_\two ^{i-1})}{3(\tau_i+\tau_{i-1})}\big(-(r_i+1)\tau_i^2+\tau_{i-1}\tau_i(r_i+1)^2\big)\\
		&=\frac{\tau_i}{3}\big(\piff  f(t_i,v_\two ^i) - \piff  f(t_{i-1},v_\two ^{i-1})\big).
	\end{align*}
	The perturbed error of the term $\mathcal{D}_{2,3}\,f(t_i,v_\two ^i)$ can be formulated by 
	\begin{align}\label{eq: D23 integral form}
		&\mathcal{D}_{2,3}\,f(t_i,v_\two ^i)-\mathcal{D}_{2,3}\,f(t_i,\bar{v}_\two ^i)\nonumber\\
		&=\frac{\tau_i}{3}\big(\piff  f(t_i,v_\two ^i)-\piff  f(t_i,\bar{v}_\two ^i) - \piff  f(t_{i-1},v_\two ^{i-1})+ \piff  f(t_{i-1},\bar{v}_\two ^{i-1})\big)\nonumber\\
		&=\frac{\tau_i}{3}\bigg(-\int_0^1f'\big(\tilde{q}_\two ^i\big)\partial_\tau\tilde{v}_\two ^{i}\zd s_1-\int_0^1\int_0^1f''\big(\tilde{g}_\two ^i\big)\partial_\tau\tilde{q}_\two ^i\tilde{v}_\two ^{i-1}\zd s_1\zd s_2\nonumber\\
		&\quad+\int_0^1f'\big(\tilde{q}_\two ^{i-1}\big)\partial_\tau\tilde{v}_\two ^{i-1}\zd s_1+\int_0^1\int_0^1f''\big(\tilde{g}_\two ^{i-1}\big)\partial_\tau\tilde{q}_\two ^{i-1}\tilde{v}_\two ^{i-2}\zd s_1\zd s_2\bigg),
	\end{align}
	where $\tilde{q}_\two ^i:=(1-s_1)v_\two ^i+s_1\bar{v}_\two ^{i}$ and $\tilde{g}_\two ^i:=(1-s_2)\tilde{q}_\two ^{i-1}+s_2\tilde{q}_\two ^{i}$. If the nonlinear function 
	$f\in C^2$, using the solution of BDF2 scheme is bounded, it is easy to find that the perturbed error  $\mathcal{D}_{2,3}\,f(t_i,v_\two ^i)-\mathcal{D}_{2,3}\,f(t_i,\bar{v}_\two ^i)$ is bounded. Thus we have the following theorem.
	\begin{theorem}\label{thm: BDF2-DC3 stability}
		If $\tau\le 1/(4L_f)$ and $f\in C^2$, the solution of \eqref{eq: BDF2-DC3 DOC perturbed equation} satisfies
		\begin{align*}
			\abs{\tilde{v}_\three^{n}}&\leq4\exp(4L_ft_{n-1})
			\braB{\abs{\tilde{v}_\three^{1}}+\tau\abs{\piff \tilde{v}_\three^1}
				+t_n\max_{2\leq i\leq n}\Big\{\abs{\varepsilon_\three^i}
				+\absb{\mathcal{D}_{2,3}f(t_i,v_\two ^i)-\mathcal{D}_{2,3}f(t_i,\bar{v}_\two ^i)}\Big\}},
		\end{align*}
		for $2\le n\le N$. Thus the variable-step BDF2-DC3 scheme \eqref{eq: BDF2-DC3 scheme} is mesh-robustly stable.
	\end{theorem}
	\begin{proof} 
		Multiplying both sides of \eqref{eq: BDF2-DC3 DOC perturbed equation} with $2\tau_n\tilde{v}_\three^n$, one
		uses Lemma \ref{lem: BDF2 orthogonal formula} to obtain that
		\begin{align*}
			\abs{\tilde{v}_\three^n}^2-\abs{\tilde{v}_\three^{n-1}}^2&+\big|\tau_n\piff  \tilde{v}_\three^n\big|^2
			\le2\tau_n \abs{\tilde{v}_\three^n}\big|\mathcal{I}_{1}^n[\tilde{v}_\three]\big|+2\tau_n\abs{\tilde{v}_\three^n}\sum_{i=2}^n\vartheta_{n-i}^{(n)}\absb{f(t_i,\bar{v}_\three^i)-f(t_i,v_\three^i)}\\
			&\quad+2\tau_n\abs{\tilde{v}_\three^n}\max_{2\leq i\leq n}\Big\{\absb{\mathcal{D}_{2,3}f(t_i,v_\two ^i)-\mathcal{D}_{2,3}f(t_i,\bar{v}_\two ^i)}+\abs{\varepsilon_\three^i}\Big\}\sum_{i=2}^n\vartheta_{n-i}^{(n)}\\
			&\le2\tau_n \abs{\tilde{v}_\three^n}\big|\mathcal{I}_{1}^n[\tilde{v}_\three]\big|+2L_f\tau_n\abs{\tilde{v}_\three^n}\sum_{i=2}^n\vartheta_{n-i}^{(n)}\abs{ \tilde{v}_\three^i}\\
			&\quad+2\tau_n\abs{\tilde{v}_\three^n}\max_{2\leq i\leq n}\Big\{\absb{\mathcal{D}_{2,3}f(t_i,v_\two ^i)-\mathcal{D}_{2,3}f(t_i,\bar{v}_\two ^i)}
			+\abs{\varepsilon_\three^i}\Big\}.
		\end{align*}
	By dropping the nonnegative term at the left side and summing $n$ in the last inequality 
	from $n=2$ to $m$ ($2\le m\le N$), one has
		\begin{align*}
			\abs{\tilde{v}_\three^m}^2&\le\,\big|\tilde{v}_\three^{1}\big|^2
			+2\sum_{j=2}^m\tau_j \big|\tilde{v}_\three^j\big|\big|\mathcal{I}_{1}^j[\tilde{v}_\three]\big|
			+2L_f\sum_{j=2}^m\tau_j\big|\tilde{v}_\three^j\big|
			\sum_{i=2}^j\vartheta_{j-i}^{(j)}\abs{\tilde{v}_\three^i}\\
			&\quad
			+2\max_{2\leq i\leq m}\Big\{\absb{\mathcal{D}_{2,3}f(t_i,v_\two ^i)-\mathcal{D}_{2,3}f(t_i,\bar{v}_\two ^i)}
			+\abs{\varepsilon_\three^i}\Big\}\sum_{j=2}^m\tau_j\big|\tilde{v}_\three^j\big|.
		\end{align*}
	Consider a fixed $m_0$ $(1\le m_0\le m)$ such that $\abs{\tilde{v}_\three^{m_0}}:=\max_{1\leq k \leq m}\abs{\tilde{v}_\three^{k}}$. We take $m:=m_0$ in the above inequality to get
		\begin{align*}
			\absb{\tilde{v}_\three^{m_0}}^2&\le\,\absb{\tilde{v}_\three^{1}}\absb{\tilde{v}_\three^{m_0}}+2\absb{\tilde{v}_\three^{m_0}}\sum_{j=2}^{m_0}\tau_j\big|\mathcal{I}_{1}^j[\tilde{v}_\three]\big|
			+2L_f\absb{\tilde{v}_\three^{m_0}}\sum_{j=2}^{m_0}\tau_j\big|\tilde{v}_\three^j\big|\\
			&\quad+2\absb{\tilde{v}_\three^{m_0}}\max_{2\leq i\leq m_0}\Big\{\absb{\mathcal{D}_{2,3}f(t_i,v_\two ^i)-\mathcal{D}_{2,3}f(t_i,\bar{v}_\two ^i)}+\abs{\varepsilon_\three^i}\Big\}\sum_{j=2}^{m_0}\tau_j.
		\end{align*}
	Removing $\absb{\tilde{v}_\three^{m_0}}$ on both sides simultaneously, one can obtain
		\begin{align}\label{proof BDF2-DC3 stability}
			\absb{\tilde{v}_\three^{m}}\le \absb{\tilde{v}_\three^{m_0}}
			&\le\absb{\tilde{v}_\three^{1}}
			+2\sum_{j=2}^m\tau_j\big|\mathcal{I}_{1}^j[\tilde{v}_\three]\big|
			+2L_f\sum_{j=2}^{m}\tau_j\big|\tilde{v}_\three^j\big|\nonumber\\
			&\quad+2t_m\max_{2\leq i\leq m}\Big\{\absb{\mathcal{D}_{2,3}f(t_i,v_\two ^i)-\mathcal{D}_{2,3}f(t_i,\bar{v}_\two ^i)}+\abs{\varepsilon_\three^i}\Big\}.
		\end{align}
	According to Lemma \ref{lem: BDF2 orthogonal formula}, the starting effect term $\mathcal{I}_{1}^n[\tilde{v}_\three]$ 
	in \eqref{eq: initial effect BDFk-DOC} can be written as
		\begin{align}\label{proof starting error BDF2}
			\mathcal{I}_{1}^n[\tilde{v}_\three]
			=\vartheta_{n-2}^{(n)}d^{(2)}_{1}\,\piff \tilde{v}_\three^1
			=-\piff \tilde{v}_\three^1\,\prod_{i=2}^n\frac{r_{i}}{1+2r_{i}}
		\end{align}
	such that
		\begin{align}\label{proof starting error DC3}
			\sum_{j=2}^m\tau_j\absb{\mathcal{I}_{1}^j[\tilde{v}_\three]}
			\le\tau\abs{\piff \tilde{v}_\three^1}\sum_{j=2}^m\frac{1}{2^{j-1}}
			\le\tau\abs{\piff \tilde{v}_\three^1}\quad\text{for $2\leq m\leq N$.}
		\end{align}
	Inserting the estimate \eqref{proof starting error DC3} into \eqref{proof BDF2-DC3 stability} and recalling the step setting $\tau\le 1/(4L_f)$, we have
		\begin{align*}
			\abs{\tilde{v}_\three^{m}}&\le 2\abs{\tilde{v}_\three^{1}}+4\tau\abs{\partial_\tau \tilde{v}_\three^1}
			+4L_f\sum_{j=2}^{m-1}\tau_j\absb{\tilde{v}_\three^j}\\
			&\quad+4t_m\max_{2\leq i\leq m}\Big\{\absb{\mathcal{D}_{2,3}f(t_i,v_\two ^i)-\mathcal{D}_{2,3}f(t_i,\bar{v}_\two ^i)}+\abs{\varepsilon_\three^i}\Big\}.
		\end{align*}
	The standard discrete Gr\"{o}nwall inequality\cite[Lemma 3.1]{LiaoZhang:2021} leads to
		\begin{align*}
			\abs{\tilde{v}_\three^{n}}&\leq4\exp(4L_ft_{n-1})\braB{\abs{\tilde{v}_\three^{1}}+\tau\abs{\piff \tilde{v}_\three^1}+t_n\max_{2\leq i\leq n}\Big\{\abs{\varepsilon_\three^i}
			+\absb{\mathcal{D}_{2,3}f(t_i,v_\two ^i)-\mathcal{D}_{2,3}f(t_i,\bar{v}_\two ^i)}\Big\}}
		\end{align*}
	for $2\leq n\leq N$. It completes the proof.
	\end{proof}

	\subsection{Convergence of the BDF2-DC3 scheme}
	This subsection proves the convergence of BDF2-DC3 scheme \eqref{eq: BDF2-DC3 scheme}. Set numerical errors 
	$$e_{\rmq}^0:=0,\quad e_{\rmq}^i:=v(t_i)-v_{\rmq}^i\quad \text{for $\rmq=\mathbf{2},\mathbf{3}$ and $i=1,\cdots,N.$}$$
	Subtracting \eqref{con: ODE DC equation} from \eqref{eq: DC2 scheme}, one gets the following error equation of the BDF2 solution
		\begin{align}\label{eq: BDF2 error equation}
			D_2 e_\two ^{i} = f(t_i,v(t_i))-f(t_i,v_\two ^i)+R_\two ^i\quad \text{for $i\geq2$,}
		\end{align}
	where the truncation error $R_\two ^i$ reads
		\begin{align}
			R_\two ^i:=-\frac{1}{6}\tau_i(\tau_i+\tau_{i-1})v'''(t_i)+O(\tau^3).
		\end{align}
	Then Lemma \ref{lem: BDF2 stability} implies the following estimate, cf. \cite{KangLiaoWang:2023cnsns,LiaoJiWangZhang:2021}. 
		\begin{theorem}\label{thm: error estimation of BDF2}
			Let the solution $v\in C^3((0,T])$. If $\tau\le 1/(4L_f)$, the  solution $v_\two ^n$ of the BDF2 scheme \eqref{eq: DC2 scheme} is mesh-robustly convergent in the sense that 
			\begin{align*} 
				\abs{e_\two ^n}\leq&\,4\exp(4L_ft_{n-1})
				\braB{\abs{e_\two ^{1}}+\tau\absb{\piff e_\two ^1}+Ct_n\tau^2}
				\quad\text{for $2\le n\leq N$.}
			\end{align*}
		\end{theorem}
		\begin{remark}\label{remark: BDF2 starting}
			As is well-known, to obtain a second-order accurate solution, $e_\two ^{1}=O(\tau^2)$, and a first-order accurate time difference quotient, $\piff e_\two ^1=O(\tau)$, some first-order starting scheme such as the BDF1 scheme is enough  \cite{KangLiaoWang:2023cnsns,LiaoJiWangZhang:2021}, especially when the problem does not require a very small starting step $\tau_1$.  In practical applications, the initial time step $\tau_1$ is typically taken to be small in order to minimize the effect of the first-order error. However, the solution computed in this way generates a very large error for stiff problems. As shown in \cite{Nishikawa2019}, the
			variable-step BDF2 method combined with any first- or higher-order one-step method is asymptotically
			equivalent to the trapezoidal method and therefore not L-stable. Thus the BDF2 method \eqref{eq: DC2 scheme} is suggested to start by a second-order L-stable RK method, which generates a third-order accurate solution at the first time level, that is, $\abs{e_\two ^{1}}+\tau\absb{\piff e_\two ^1}=O(\tau^3)$. 
		\end{remark}
	To accurately describe the effect of initial errors, we introduce some time-level-dependent quantities $\sigma_{\mathbf{q}}^n$  $(\mathbf{q}=\mathbf{2},\mathbf{3},\mathbf{4})$, 
	which are asymptotically vanished as the index $n$ is properly large, 
		\begin{align} \label{def: decaying factors}
			\sigma_\two ^n:=&\,\prod_{\ell=2}^n\frac{r_{\ell}}{1+2r_{\ell}},\quad
			\sigma_\three^n:=\sum\limits_{i=2}^n\vartheta_{n-i}^{(n)}\sigma_\two ^i,\quad
			\sigma_\four^n:=\sum\limits_{i=2}^n\vartheta_{n-i}^{(n)}\sigma_\three^i.
		\end{align}
	One can obtain the following upper bounds
		\begin{align} \label{ieq: bound decaying factors}
			\sigma_\two ^n\le 2^{1-n},\quad
			\sigma_\three^n\le 2^{1-n}n\quad\text{and}\quad
			\sigma_\four^n\le 2^{-n}n^2.
		\end{align}
	Actually, by using Lemma \ref{lem: BDF2 orthogonal formula}, it is not difficult to find that
		\begin{align*}
			\sigma_\three^n=&\,
			\prod_{\ell=2}^n\frac{r_{\ell}}{1+2r_{\ell}}\sum\limits_{i=2}^n\frac{1}{d_0^{(i)}}
			=\sigma_\two ^n\sum\limits_{i=2}^n\frac{1}{d_0^{(i)}}\le 2^{1-n}(n-1),\\
			\sigma_\four^n=&\,\sum\limits_{i=2}^n\vartheta_{n-i}^{(n)}
			\sum\limits_{k=2}^i\vartheta_{i-k}^{(i)}\sigma_\two ^k
			=\sum\limits_{k=2}^n\sigma_\two ^k\sum\limits_{i=k}^n\vartheta_{n-i}^{(n)}\vartheta_{i-k}^{(i)}
			=\sum\limits_{k=2}^n\sigma_\two ^k
			\sum\limits_{i=k}^n\frac{1}{d_0^{(i)}d_0^{(k)}}\prod_{\ell=k+1}^n\frac{r_{\ell}}{1+2r_{\ell}}\\
			=&\,\prod_{\ell=2}^n\frac{r_{\ell}}{1+2r_{\ell}}\sum\limits_{k=2}^n\sum\limits_{i=k}^n\frac{1}{d_0^{(i)}d_0^{(k)}}
			\leq 2^{-n}(n-1)n.
		\end{align*}

	Since the BDF2 solution and its  time difference quotient are involved 
	in the variable-step BDF2-DC3 scheme \eqref{eq: BDF2-DC3 scheme}, 
	we present the following result on the discrete time derivative.
	\begin{lemma}\label{lemma: BDF2 time derivative error}
	Let the solution $v\in C^3((0,T])$ and the nonlinear function $f\in C^2$. If $\tau\le 1/(4L_f)$, the solution $v_\two ^n$ of the BDF2 scheme \eqref{eq: DC2 scheme} satisfy
		\begin{align*} 	
			\abs{\piff e_\two ^n}
			\leq&\,\sigma_\two ^n\abs{\piff e_\two ^1}
			+C\bra{\abs{e_\two ^{1}}+\tau\absb{\piff e_\two ^1}+t_n\tau^2}
			\quad\text{for $2\le n\leq N$.}
		\end{align*}
	\end{lemma}
	\begin{proof}		
		Multiplying both sides of the error equation \eqref{eq: BDF2 error equation} by the DOC kernels
		$\vartheta_{n-i}^{(n)}$, and summing $i$ from 2 to $n$, we apply \eqref{eq: BDF2-DOC transform} to get
		\begin{align*}
			\piff e_\two ^n
			=-\mathcal{I}_{1}^n[e_\two ]+\sum\limits_{i=2}^n\vartheta_{n-i}^{(n)}
			\big\{f(t_i,v(t_i))-f(t_i,v_\two ^i)+R_\two ^i\big\}\quad \text{for $n\geq 2$}.
		\end{align*}
	Recalling the definition \eqref{eq: initial effect BDFk-DOC}, one applies Lemma \ref{lem: BDF2 orthogonal formula} and Theorem \ref{thm: error estimation of BDF2} to obtain
		\begin{align*}
			\abs{\piff e_\two ^n}
			\leq&\,\sigma_\two ^n\abs{\piff e_\two ^1}+\max_{2\leq i\leq n}\big\{L_f\abs{e_\two ^i}+\abs{R_\two ^i}\big\}\\
			\le&\,\sigma_\two ^n\abs{\piff e_\two ^1}+4L_f\exp(4L_ft_{n-1})
			\braB{\abs{e_\two ^{1}}+\tau\absb{\piff e_\two ^1}+Ct_n\tau^2}
		\end{align*}
	for $2\le n\leq N$. It completes the proof.		
	\end{proof} 
	Lemma \ref{lemma: BDF2 time derivative error} says that the time difference quotient is asymptotically second-order accurate (as the level index $n$ is properly large) if the BDF1 method or the second-order RK method is used to compute the first-level solution $v_\two ^{1}$ in the sense that, cf. the data in Table \ref{table: partial error graded meshes for example1},
	 \begin{align*}
	 	\abs{\piff e_\two ^n}\le
	 	C\bra{\abs{e_\two ^{1}}+\brat{\tau+2^{-n}}\absb{\piff e_\two ^1}+t_n\tau^2}
	 	\quad\text{for $2\le n\leq N$.}
	 \end{align*}

	As an updated version of \cite[Lemma 3.1]{BourgaultGaron:2022}, the following lemma considers the error effect of the third-order correction term $\mathcal{D}_{2,3} f(t_n,v_\two ^n)$.
	\begin{lemma}\label{lem: truncation error of DC3}
		Assume the solution $v(t)\in C^3((0,T])$ and the nonlinear function $f\in C^2$. If $\tau\le 1/(4L_f)$, it holds that
		\begin{align*} 
			\sum\limits_{i=2}^n\vartheta_{n-i}^{(n)}\absb{\mathcal{D}_{2,3} f(t_i,v_\two ^i)- \mathcal{D}_{2,3} f(t_i,v(t_i))}\leq C\tau\braB{\abs{e_\two ^{1}}
				+(\tau+\sigma_\three^n)\abs{\piff e_\two ^1}+t_n\tau^2},\quad2\le n\leq N.
		\end{align*}
	\end{lemma}
	\begin{proof}	
		Recalling the integral form \eqref{eq: D23 integral form}, one has
		\begin{align*}
			\mathcal{D}_{2,3}\,f(t_i,v_\two ^i)&-\mathcal{D}_{2,3}\,f(t_i,v(t_i))\nonumber\\
			&=\frac{\tau_i}{3}\bigg(\int_0^1f'\big(q_\two ^i\big)\partial_\tau e_\two ^{i}\zd s_1+\int_0^1\int_0^1f''\big(g_\two ^i\big)\partial_\tau q_\two ^ie_\two ^{i-1}\zd s_1\zd s_2\nonumber\\
			&\quad-\int_0^1f'\big(q_\two ^{i-1}\big)\partial_\tau e_\two ^{i-1}\zd s_1-\int_0^1\int_0^1f''\big(g_\two ^{i-1}\big)\partial_\tau q_\two ^{i-1}e_\two ^{i-2}\zd s_1\zd s_2\bigg),
		\end{align*}
		where $q_\two ^i:=(1-s_1)v(t_i)+s_1v_\two ^i$ and $g_\two ^i:=(1-s_2)q_\two ^{i-1}+s_2q_\two ^{i}$. With the help of Theorem \ref{thm: error estimation of BDF2} and Lemma \ref{lemma: BDF2 time derivative error}, it follows that
		\begin{align*}
			&\absb{\mathcal{D}_{2,3}\,f(t_i,v_\two ^i)-\mathcal{D}_{2,3}\,f(t_i,v(t_i))}\\
			&\leq  C\tau_i\braB{\abs{\partial_\tau e_\two ^{i}}+\abs{\partial_\tau e_\two ^{i-1}}
			+\abs{e_\two ^{i-1}}+\abs{e_\two ^{i-2}}}\\
		 &\leq C\tau\sigma_\two ^i\abs{\piff e_\two ^1}+C\tau\braB{\abs{e_\two ^{1}}+\tau\absb{\piff e_\two ^1}+t_i\tau^2}
		 +C\tau\exp(4L_ft_{i-1})
		 \braB{\abs{e_\two ^{1}}+\tau\absb{\piff e_\two ^1}+t_i\tau^2}\\
		 &\leq C\tau\sigma_\two ^i\abs{\piff e_\two ^1}
		 +C\tau\exp(4L_ft_{i-1})
		 \braB{\abs{e_\two ^{1}}+\tau\absb{\piff e_\two ^1}+t_i\tau^2}.
		\end{align*}
		Then Lemma \ref{lem: BDF2 orthogonal formula} together with the definition \eqref{def: decaying factors} gives
		\begin{align*} 
			\sum\limits_{i=2}^n\vartheta_{n-i}^{(n)}\absb{\mathcal{D}_{2,3} f(t_i,v_\two ^i)- \mathcal{D}_{2,3} f(t_i,v(t_i))}\leq C\tau\sigma_\three^n\abs{\piff e_\two ^1}+  C\braB{\tau\abs{e_\two ^{1}}+\tau^2\abs{\piff e_\two ^1}+t_n\tau^3}
		\end{align*}
		for $2\le n\leq N$. It completes the proof.
	\end{proof}

	 We are in a position to present the convergence result for the BDF2-DC3 scheme \eqref{eq: BDF2-DC3 scheme}. 
	\begin{theorem}\label{thm: error estimation of BDF2-DC3}
		Let $v\in C^4((0,T])$ and $v_\two ^n$ be the solution of BDF2 scheme \eqref{eq: DC2 scheme}. If $\tau\le 1/(4L_f)$ and $f\in C^2$, the  solution $v_\three^n$ of the BDF2-DC3 scheme \eqref{eq: BDF2-DC3 scheme} is mesh-robustly convergent, 
		\begin{align*} 
			\abs{e_\three^n}\leq C\exp(4L_ft_{n-1})
			\Big(\abs{e_\three^{1}}+\tau\absb{\piff e_\three^1}
				+\tau\abs{e_\two ^1}+\tau(\tau+\sigma_\three^n)\abs{\piff e_\two ^1}+t_n\tau^3\Big),		\quad\text{$2\le n\leq N$.}		
		\end{align*}	
	\end{theorem}
	\begin{proof}
		By 	subtracting \eqref{con: ODE DC equation} from \eqref{eq: BDF2-DC3 scheme}, we have the  following error equation of \eqref{eq: BDF2-DC3 scheme},
		\begin{align}\label{eq: BDF2-DC3 error equation}
			D_2 e_\mathbf{3}^{i} = f(t_i,v(t_i))-f(t_i,v_\mathbf{3}^i)+\mathcal{D}_{2,3} f(t_i,v_\mathbf{2}^i)- \mathcal{D}_{2,3} f(t_i,v(t_i))+R_\mathbf{3}^i \quad\text{for $i\geq2$,}
		\end{align}
		where the truncation error $R_\three^i$ reads
		\begin{align*}
			R_\three^i:=&\,\mathcal{D}_{2,3} f(t_i,v(t_i))-R^i(t_i)
			=-\frac{1}{4!}\big(r_i(\tau_i+\tau_{i-1})^3-(1+r_i)\tau_i^3\big)v^{(4)}(t_i)+O(\tau^4)\nonumber\\
			=&\,-\frac{1}{24}\tau_i(2\tau_i^2+3\tau_i\tau_{i-1}+\tau_{i-1}^2)v^{(4)}(t_i)+O(\tau^4).
		\end{align*}
		Multiplying both sides of \eqref{eq: BDF2-DC3 error equation} by the DOC kernels
		$\vartheta_{n-i}^{(n)}$, and summing $i$ from 2 to $n$, one gets
		\begin{align}\label{eq: BDF2-DC3 DOC error equation}
			\piff e_\three^n&
			=-\mathcal{I}_{1}^n[e_\three]+\sum\limits_{i=2}^n\vartheta_{n-i}^{(n)}
			\braB{f(t_i,v(t_i))-f(t_i,v_\three^i)+R_\three^i}\nonumber\\
			&\quad+\sum\limits_{i=2}^n\vartheta_{n-i}^{(n)}\braB{\mathcal{D}_{2,3}f(t_i,v_\two ^i)-\mathcal{D}_{2,3}f(t_i,v(t_i))}\quad \text{for $n\geq 2$}.
		\end{align}	
		Then, with the help of Lemma \ref{lem: truncation error of DC3}, it is not difficult to achieve the desired error estimate by following the proof of Theorem \ref{thm: BDF2-DC3 stability}.
	\end{proof}
	Theorem \ref{thm: error estimation of BDF2-DC3} together with the upper bounds in \eqref{ieq: bound decaying factors} suggests that the BDF2-DC3 scheme is mesh-robustly third-order convergent if the starting solutions are accurate enough such that
	\begin{align*} 
	\abs{e_\three^{1}}+\tau\absb{\piff e_\three^1}
			+\tau\abs{e_\two ^1}+\tau\brat{\tau+2^{-n}n}\abs{\piff e_\two ^1}=O(\tau^3).		
	\end{align*}
	\begin{remark}\label{remark: BDF2-DC3 starting}
		According to Remark \ref{remark: BDF2 starting}, one has $\abs{e_\two ^1}+\tau\abs{\piff e_\two ^1}=O(\tau^3)$ 
		by using the second-order L-stable RK method to start the BDF2 scheme. Similarly,  we have
		$\abs{e_\three^1}+\tau\abs{\piff e_\three^1}=O(\tau^3)$ if 
		the second-order RK method is used to compute the starting solution $v_\three^1$ for the BDF2-DC3 scheme 
		\eqref{eq: BDF2-DC3 scheme}. As seen, our error estimate in Theorem \ref{thm: error estimation of BDF2-DC3} reflects the numerical effect on the solutions for different choices of the starting methods.  We list three different cases for different starting groups (here and hereafter,  the tuple ``\textbf{Scheme A+Scheme B}" represents that the starting solutions $v_\two ^1$ and $v_\three^1$ are generated by \textbf{Scheme A} and \textbf{Scheme B}, respectively):
		\begin{description}
			\item[(RK2+RK2)] The BDF2 solution is second-order accurate and the BDF2-DC3 scheme  \eqref{eq: BDF2-DC3 scheme} generates the third-order solution.
			\item[(BDF1+RK2)] For the problems that do not require a very small starting step $\tau_1$, one can use the BDF1 scheme and the second-order L-stable RK method to compute the starting solutions $v_\two ^1$ and $v_\three^1$, respectively. In such case, Theorem \ref{thm: error estimation of BDF2-DC3} suggests that the variable-step BDF2-DC3 scheme is also third-order convergent on arbitrary time meshes because  $\tau\abs{e_\two ^1}+\tau\brat{\tau+2^{-n}n}\abs{\piff e_\two ^1}$ is asymptotically third-order accurate. 
			\item[(BDF1+BDF1)] If $v_\two ^1$ and $v_\three^1$ are computed by the  BDF1 scheme, Theorem \ref{thm: error estimation of BDF2-DC3} suggests that the BDF2-DC3 solution is only second-order accurate since  $\abs{e_\three^1}+\tau\abs{\piff e_\three^1}=O(\tau^2)$. 
		\end{description}
		It seems that the BDF2-DC3 method has no aftereffect, which means that the first-order starting method of BDF2 scheme will not cause a loss in the accuracy of the BDF2-DC3 scheme. That is to say, when the available BDF2 solutions are  used to obtain the third-order BDF2-DC3 solution, there is no need to recompute the BDF2 starting value by certain second-order method.	
	\end{remark}

	\section{Stability and convergence of BDF2-DC3-DC4 scheme}
	\setcounter{equation}{0}
	\subsection{ Stability of the BDF2-DC3-DC4 scheme}
	The stability of BDF2-DC3-DC4 scheme \eqref{eq: BDF2-DC3-DC4 scheme} is also studied via its convolutional form.
	Multiplying both sides of the perturbed error equation \eqref{eq: BDF2-DC3-DC4 perturbed equation} by the DOC kernels
	$\vartheta_{n-i}^{(n)}$, and summing $i$ from 3 to $n$, we apply the identity \eqref{eq: BDF2-DOC identity} to get
	\begin{align}\label{eq: BDF2-DC3-DC4 DOC error equation}
		\piff \tilde{v}_\four^n
		=-\mathcal{I}_{2}^n[\tilde{v}_\four]+\sum\limits_{i=3}^n\vartheta_{n-i}^{(n)}
		\kbraB{f(t_i,\bar{v}_\four^i)-f(t_i,v_\four^i)+\mathcal{D}_{2,4}f(t_i,v_\three^i)-\mathcal{D}_{2,4}f(t_i,\bar{v}_\three^i)+\varepsilon_\four^i}
	\end{align}
	for $3\leq n\leq N$, where $\mathcal{I}_{2}^n[\tilde{v}_\four]$ represents the starting effect on the BDF2-DC3-DC4 solution error at the time $t_n$ (reminding the fact $d^{(i)}_{\ell}=0$ for $2\le \ell\le i$)
	\begin{align}\label{eq: initial effect BDF2-DC3-DC4}
		\mathcal{I}_{2}^n[\tilde{v}_\four]
		:=&\,\sum_{i=3}^n\vartheta_{n-i}^{(n)}
		\sum_{j=1}^2d^{(i)}_{i-j}\piff \tilde{v}_\four^j
		=\vartheta_{n-3}^{(n)}d^{(3)}_{1}\piff \tilde{v}_\four^2
		\quad\text{for $3\le n\le N$.}
	\end{align}
	By the definition \eqref{de: D24}, simple calculations give
	\begin{align*}
		\mathcal{D}_{2,4}\,f(t_i,v_\three^i)&=\mathcal{D}_{2,3}\,f(t_i,v_\three^i)+\mathcal{D}_{3,4}\,f(t_i,v_\three^i)\\
		&=\mathcal{D}_{2,3}\,f(t_i,v_\three^i)+\frac{\tau_i(\tau_i+\tau_{i-1})(2\tau_i+\tau_{i-1})}{12(\tau_i+\tau_{i-1}+\tau_{i-2})}\\
		&\quad\times\bigg[\frac{\piff  f(t_{i},v_\three^{i})-\piff  f(t_{i-1},v_\three^{i-1})}{t_i-t_{i-2}}-\frac{\partial_\tau f(t_{i-1},v_\three^{i-1})-\partial_\tau f(t_{i-2},v_\three^{i-2})}{t_{i-1}-t_{i-3}}\bigg].
	\end{align*} 
	The perturbed error of the term $\mathcal{D}_{2,4}\,f(t_i,v_\three^i)$ can be formulated by 
	\begin{align}\label{eq: D24 perturb error integral form}
		\absb{\mathcal{D}_{2,4}\,f(t_i,v_\three^i)-\mathcal{D}_{2,4}\,f(t_i,\bar{v}_\three^i)}
		&\leq\absb{\mathcal{D}_{2,3}\,f(t_i,v_\three^i)-\mathcal{D}_{2,3}\,f(t_i,\bar{v}_\three^i)}
		+\frac{\tau_i}{6}\big|\tilde{F}^{i}_\three+\tilde{K}^{i}_\three
		- \tilde{F}^{i-1}_\three-\tilde{K}^{i-1}_\three\big|\nonumber\\
		&\quad
		+\frac{\tau_ir_{i-1}(r_i+1)}{6(r_{i-1}+1)}\big|\tilde{F}^{i-1}_\three+\tilde{K}^{i-1}_\three
		-\tilde{F}^{i-2}_\three-\tilde{K}^{i-2}_\three\big|,
	\end{align}
	where the simple facts, $\frac{\tau_i+\tau_{i-1}}{\tau_i+\tau_{i-1}+\tau_{i-2}}<1$ and $\frac{2\tau_i+\tau_{i-1}}{12}<\frac{\tau_i+\tau_{i-1}}{6}$, have been used, and
	\begin{align*}
		\tilde{F}^{i}_\three:=&\,\int_0^1f'\big(\tilde{q}_\three^{i}\big)\partial_\tau\tilde{v}_\three^i\zd s_1,\qquad\qquad\qquad\,\,\, \tilde{q}_\three^i:=(1-s_1)v_\three^i+s_1\bar{v}_\three^i,\\
		\tilde{K}^{i}_\three:=&\,\int_0^1\int_0^1f''\big(\tilde{g}_\three^{i}\big)\piff \tilde{q}_\three^{i}\tilde{v}_\three^{i-1}\zd s_1\zd s_2,\quad\tilde{g}_\three^i:=(1-s_2)\tilde{q}_\three^{i-1}+s_2\tilde{q}_\three^{i}.
	\end{align*} 
	Since the solution of BDF2-DC3 scheme and the perturbed error \eqref{eq: D23 integral form} of the term $\mathcal{D}_{2,3}\,f(t_i,v_\two ^i)$ are bounded, it is not difficult to check that the perturbed error $\absb{\mathcal{D}_{2,4}\,f(t_i,v_\three^i)-\mathcal{D}_{2,4}\,f(t_i,\bar{v}_\three^i)}$ is also bounded. Thus we have the following theorem.
	\begin{theorem}\label{thm: BDF2-DC3-DC4 stability}
		If $\tau\le 1/(4L_f)$ and $f\in C^2$, the solution of \eqref{eq: BDF2-DC3-DC4 DOC error equation} satisfies
		\begin{align*}
			\abs{\tilde{v}_\four^{n}}\leq4\exp(4L_ft_{n-1})
			\Big(\abs{\tilde{v}_\four^{2}}
				+\tau\abs{\piff \tilde{v}_\four^2}
				+t_n\max_{3\leq i\leq n}\Big\{\abs{\varepsilon_\four^i}
				+\absb{\mathcal{D}_{2,4}f(t_i,v_\three^i)-\mathcal{D}_{2,4}f(t_i,\bar{v}_\three^i)}\Big\}\Big)
		\end{align*}
	for $3\le n\le N$. Thus the BDF2-DC3-DC4 scheme \eqref{eq: BDF2-DC3-DC4 scheme} is mesh-robustly stable.
	\end{theorem}
	\begin{proof} 
	The proof is similar to Theorem \ref{thm: BDF2-DC3 stability}, so it will not be repeated here.
	\end{proof}

	\subsection{Convergence of the BDF2-DC3-DC4 scheme}
	Since the BDF2-DC3 solution and its time difference quotient are involved in the variable-step BDF2-DC3-DC4 scheme \eqref{eq: BDF2-DC3-DC4 scheme}, we need the following result on the discrete time derivative.
	\begin{lemma}\label{lemma: BDF2-DC3 time derivative error}
		Assume the solution $v(t)\in C^4((0,T])$ and the nonlinear function $f\in C^2$. If $\tau\le 1/(4L_f)$, the discrete solution $v_\three^n$ of the variable-step BDF2-DC3 scheme \eqref{eq: BDF2-DC3 scheme} satisfies
		\begin{align*} 	
			\abs{\piff e_\three^n}
			\leq&\,C\Big(\abs{e_\three^{1}}+\brat{\tau+\sigma_\two ^n}\absb{\piff e_\three^1}
			+\tau\abs{e_\two ^1}+\tau\brat{\tau+\sigma_\three^n}\abs{\piff e_\two ^1}+t_n\tau^3\Big)
			\quad\text{for $2\le n\leq N$.}
		\end{align*}
	\end{lemma}
	\begin{proof}Reminding the definition \eqref{eq: initial effect BDFk-DOC} and 
		Lemma \ref{lem: BDF2 orthogonal formula}, one derives from 
		the error equation \eqref{eq: BDF2-DC3 DOC error equation} of the BDF2-DC3 scheme that
		\begin{align*}
			\abs{\piff e_\three^n}
			\le&\, \absb{\mathcal{I}_{1}^n[e_\three]}
			+\sum\limits_{i=2}^n\vartheta_{n-i}^{(n)}
			\absb{\mathcal{D}_{2,3}f(t_i,v_\two ^i)-\mathcal{D}_{2,3}f(t_i,v(t_i))}
			+\max_{2\leq i\leq n}\Big\{L_f\abs{e_\three^i}+\abs{R_\three^i}\Big\}			\\
			\le&\, \sigma_\two ^n\abs{\piff e_\three^1}
			+C\tau\braB{\abs{e_\two ^{1}}
				+(\tau+\sigma_\three^n)\abs{\piff e_\two ^1}+t_n\tau^2}\nonumber\\
			&\,	+C\exp(4L_ft_{n-1})
			\Big(\abs{e_\three^{1}}+\tau\absb{\piff e_\three^1}
			+\tau\abs{e_\two ^1}+\tau(\tau+\sigma_\three^n)\abs{\piff e_\two ^1}+t_n\tau^3\Big)\\
			\le&\, \sigma_\two ^n\abs{\piff e_\three^1}
			+C\Big(\abs{e_\three^{1}}+\tau\absb{\piff e_\three^1}
			+\tau\abs{e_\two ^1}+\tau\brat{\tau+\sigma_\three^n}\abs{\piff e_\two ^1}+t_n\tau^3\Big),
		\end{align*}
		where Lemma \ref{lem: truncation error of DC3} and Theorem \ref{thm: error estimation of BDF2-DC3} have been used. It completes the proof.		
	\end{proof} 
	\begin{remark}\label{remark: error time difference quotient}
		Lemma \ref{lemma: BDF2-DC3 time derivative error} says that the time difference quotient $\piff e_\three^n$ is asymptotically third-order accurate (as the level index $n$ is properly large) if the RK2 method is used to compute the first-level solution $v_\three^{1}$ in the sense that
		\begin{align*}
			\abs{\piff e_\three^n}
			\leq&\,C\Big(\abs{e_\three^{1}}+\brat{\tau+2^{-n}}\absb{\piff e_\three^1}
			+\tau\abs{e_\two ^1}+\tau\brat{\tau+2^{-n}n}\abs{\piff e_\two ^1}+t_n\tau^3\Big).
		\end{align*}
		To illustrate the asymptotic behaviors of the discrete derivatives $\piff e_\three^n$ and $\piff e_\two ^n$,
		cf. Lemma \ref{lemma: BDF2 time derivative error}, we run the two schemes for a linear model $v'(t)=v\cos t$, see Example \ref{example1} in Section 5, on the graded mesh $t_k=(k/n)^2$. 
		In Table \ref{table: partial error graded meshes for example1}, we use the BDF1 and RK2 schemes to startup the BDF2 and BDF2-DC3 methods, respectively. As expected,  $\piff e_\two ^n$ and $\piff e_\three^n$ are asymptotically second-order and third-order accurate, respectively.
		\begin{table}[htb!]
			\begin{center} 
				\caption{Errors of discrete time derivatives of BDF2 and BDF2-DC3 methods
					\label{table: partial error graded meshes for example1}}\vspace*{0.3pt}
				\def\temptablewidth{1\textwidth}
				{\rule{\temptablewidth}{0.5pt}}
				\begin{tabular*}{\temptablewidth}{@{\extracolsep{\fill}}cccccc}
					&\multirow{2}{*}{$n$} &\multicolumn{2}{c}{BDF2} &\multicolumn{2}{c}{BDF2-DC3}   \\
					\cline{3-4}          \cline{5-6}          
					&      &$\partial_\tau e_\two ^n$     &Order  
					&$\partial_\tau e_\three^n$  &Order  \\
					\midrule         
					&100     &6.55e-04 &-    &2.27e-06 &- \\
					&200     &1.66e-04 &1.99 &2.85e-07 &3.00 \\	
					&400     &4.18e-05 &1.99 &3.58e-08 &3.00 \\
					&800     &1.05e-05 &2.00 &4.48e-09 &3.00 \\
				\end{tabular*}
				{\rule{\temptablewidth}{0.5pt}}
			\end{center}
		\end{table}	
	\end{remark}
	Now we consider the error effect of the fourth-order correction term $\mathcal{D}_{2,4} f(t_n,v_\three^n)$.
	\begin{lemma}\label{lem: truncation error of DC4}
		Assume the solution $v(t)\in C^4((0,T])$ and the nonlinear function $f\in C^2$. If $\tau\le 1/(4L_f)$, it holds that
		\begin{align*}
			&\sum\limits_{i=3}^n\vartheta_{n-i}^{(n)}\absb{\mathcal{D}_{2,4} f(t_i,v_\three^i)- \mathcal{D}_{2,4} f(t_i,v(t_i))}\\&\hspace{1cm}\leq Cr_{\max}\tau\Big(\abs{e_\three^{1}}+\brat{\tau+\sigma_\three^n}\absb{\piff e_\three^1}
			+\tau\abs{e_\two ^1}+\tau\brat{\tau+\sigma_\four^n}\abs{\piff e_\two ^1}+t_n\tau^3\Big),\quad 3\le n\leq N.
		\end{align*}	
	\end{lemma}
	\begin{proof}		
		By the integral form \eqref{eq: D24 perturb error integral form}, one obtains
		\begin{align}\label{eq: D24 error integral form}
			&\absb{\mathcal{D}_{2,4}\,f(t_i,v_\three^i)-\mathcal{D}_{2,4}\,f(t_i,v(t_i))}\nonumber\\
			&\hspace{2cm}\leq\absb{\mathcal{D}_{2,3}\,f(t_i,v_\three^i)-\mathcal{D}_{2,3}\,f(t_i,v(t_i))}
			+\frac{\tau_i}{6}\big|F^{i}_\three+K^{i}_\three
			- F^{i-1}_\three-K^{i-1}_\three\big|\nonumber\\
			&\hspace{2.5cm}+\frac{\tau_ir_{i-1}(r_i+1)}{6(r_{i-1}+1)}\big|F^{i-1}_\three+K^{i-1}_\three
			-F^{i-2}_\three-K^{i-2}_\three\big|,
		\end{align}
		where the fact, $\frac{\tau_i+\tau_{i-1}}{\tau_i+\tau_{i-1}+\tau_{i-2}}<1$ and $\frac{2\tau_i+\tau_{i-1}}{12}<\frac{\tau_i+\tau_{i-1}}{6}$, were used in the last inequlity, and
		\begin{align*} 
			F^{i}_\three&:=\int_0^1f'\big(q_\three^{i}\big)\partial_\tau e_\three^i\zd s_1,\qquad\qquad\qquad\,\,\, q_\three^i:=(1-s_1)v(t_i)+s_1v_\three^i,\\ K^{i}_\three&:=\int_0^1\int_0^1f''\big(g_\three^{i}\big)\piff q_\three^{i}e_\three^{i-1}\zd s_1\zd s_2,\quad\, g_\three^i:=(1-s_2)q_\three^{i-1}+s_2q_\three^{i}.
		\end{align*}
		Thus one has
		\begin{align*}
			\abs{\mathcal{D}_{2,4}\,f(t_i,v_\three^i)-\mathcal{D}_{2,4}\,f(t_i,v(t_i))}\le&\, \absb{\mathcal{D}_{2,3}\,f(t_i,v_\three^i)-\mathcal{D}_{2,3}\,f(t_i,v(t_i))}\\
			&\,+Cr_{\max}\tau_i\sum_{\ell=i-2}^i\absb{\partial_\tau e_\three^{\ell}}+Cr_{\max}\tau_i\sum_{\ell=i-3}^{i-1}\absb{e_\three^{\ell}},\quad i\ge3.
		\end{align*}
		By following the proof of Lemma \ref{lem: truncation error of DC3}, one has
		\begin{align*}
			&\absb{\mathcal{D}_{2,3}\,f(t_i,v_\three^i)-\mathcal{D}_{2,3}\,f(t_i,v(t_i))}		
			\leq  C\tau_i\braB{\abs{\partial_\tau e_\three^{i}}+\abs{\partial_\tau e_\three^{i-1}}
				+\abs{e_\three^{i-1}}+\abs{e_\three^{i-2}}},\quad i\ge3.
		\end{align*}
		Thus it follows from Theorem \ref{thm: error estimation of BDF2-DC3} 
		and Lemma \ref{lemma: BDF2-DC3 time derivative error} that
		\begin{align*}
			&\absb{\mathcal{D}_{2,4}\,f(t_i,v_\three^i)-\mathcal{D}_{2,4}\,f(t_i,v(t_i))}\\		
			&\leq  Cr_{\max}\tau_i\sum_{\ell=i-2}^i\absb{\partial_\tau e_\three^{\ell}}+Cr_{\max}\tau_i\sum_{\ell=i-3}^{i-1}\absb{e_\three^{\ell}}\\
			&\leq Cr_{\max}\tau\Big(\abs{e_\three^{1}}+\brat{\tau+\sigma_\two ^i}\absb{\piff e_\three^1}
			+\tau\abs{e_\two ^1}+\tau\brat{\tau+\sigma_\three^i}\abs{\piff e_\two ^1}+t_i\tau^3\Big)\\
			&\quad+Cr_{\max}\tau
			\Big(\abs{e_\three^{1}}+\tau\absb{\piff e_\three^1}
			+\tau\abs{e_\two ^1}+\tau(\tau+\sigma_\three^i)\abs{\piff e_\two ^1}+t_i\tau^3\Big)\\
			&\leq Cr_{\max}\tau\Big(\abs{e_\three^{1}}+\brat{\tau+\sigma_\two ^i}\absb{\piff e_\three^1}
			+\tau\abs{e_\two ^1}+\tau\brat{\tau+\sigma_\three^i}\abs{\piff e_\two ^1}+t_i\tau^3\Big).
		\end{align*}
		Recalling the definitions in  \eqref{def: decaying factors}, we apply Lemma \ref{lem: BDF2 orthogonal formula} to find that
		\begin{align*}
			&\sum\limits_{i=3}^n\vartheta_{n-i}^{(n)}\absb{\mathcal{D}_{2,4} f(t_i,v_\three^i)- \mathcal{D}_{2,4} f(t_i,v(t_i))}\\&\hspace{2cm}\leq Cr_{\max}\tau\Big(\abs{e_\three^{1}}+\brat{\tau+\sigma_\three^n}\absb{\piff e_\three^1}
			+\tau\abs{e_\two ^1}+\tau\brat{\tau+\sigma_\four^n}\abs{\piff e_\two ^1}+t_n\tau^3\Big).
		\end{align*}	
		It completes the proof.	
	\end{proof}
	
	Now we present the convergence of BDF2-DC3-DC4 scheme \eqref{eq: BDF2-DC3-DC4 scheme}. Let numerical errors 
	$$e_\four^0:=0,\quad e_\four^i:=v(t_i)-v_\four^i\quad \text{ $i=1,\cdots,N.$}$$ 
	Subtracting \eqref{con: ODE DC equation} form \eqref{eq: BDF2-DC3-DC4 scheme}, one gets the following error equation of BDF2-DC3-DC4 scheme
	\begin{align}\label{eq: BDF2-DC3-DC4 error equation}
		D_2 e_\mathbf{4}^{i} = f(t_i,v(t_i))-f(t_i,v_\mathbf{4}^i)+\mathcal{D}_{2,4} f(t_i,v_\mathbf{3}^i)- \mathcal{D}_{2,4} f(t_i,v(t_i))+R_\mathbf{4}^i\quad\text{for $i\geq 3,$}
	\end{align}
	where the truncation error $R_\four^i$ of BDF2-DC3-DC4 scheme is
	\begin{align*}
		R_\four^i&=\mathcal{D}_{2,4} f(t_i,v(t_i))-R^i(t_i)
		=\frac{v^{(5)}(t_i)}{5!}\big(r_i(\tau_i+\tau_{i-1})^4-(1+r_i)\tau_i^4\big)+O(\tau^5)\nonumber\\
		&=\frac{1}{120}\tau_i(3\tau_i^2+3\tau_i\tau_{i-1}+\tau_{i-1}^2)(\tau_i+\tau_{i-1})v^{(5)}(t_i)+O(\tau^5).
	\end{align*}
	Multiplying both sides of \eqref{eq: BDF2-DC3-DC4 error equation} by the DOC kernels
	$\vartheta_{n-i}^{(n)}$, and summing $i$ from 3 to $n$, one gets
	\begin{align*}
		\piff e_\four^n
		&=-\mathcal{I}_{2}^n[e_\four]+\sum\limits_{i=3}^n\vartheta_{n-i}^{(n)}
		\braB{f(t_i,v(t_i))-f(t_i,v_\four^i)+R_\four^i}\nonumber\\
		&\quad + \sum\limits_{i=3}^n\vartheta_{n-i}^{(n)}
		\braB{\mathcal{D}_{2,4}f(t_i,v_\three^i)-\mathcal{D}_{2,4}f(t_i,v(t_i))}\quad \text{for $n\geq 3$},
	\end{align*}
	where the starting effect term $\mathcal{I}_{2}^n[e_\four]$ is defined via \eqref{eq: initial effect BDF2-DC3-DC4}.
	Then, by using this error equation and Lemma \ref{lem: truncation error of DC4}, one can follow the proof of Theorem \ref{thm: BDF2-DC3-DC4 stability} to prove the following theorem.
	\begin{theorem}\label{thm: error estimation of BDF2-DC3-DC4}
		Let the solution $v\in C^5((0,T])$ and $v_\three^n$ be the discrete solution  of the BDF2-DC3 scheme \eqref{eq: BDF2-DC3 scheme}. If $\tau\le 1/(4L_f)$ and $f\in C^2$, the numerical solution $v_\four^n$ of the BDF2-DC3-DC4 scheme \eqref{eq: BDF2-DC3-DC4 scheme} is convergent in the sense that, 
		\begin{align*}
			\abs{e_\four^{n}}\leq Cr_{\max}&\,\exp(4L_ft_{n-1})
			\Big(\abs{e_\four^{2}}
			+\tau\abs{\piff e_\four^2}
			+t_n\tau^4\Big.\\
			&\,\Big.+\tau\abs{e_\three^{1}}
			+\tau\brat{\tau+\sigma_\three^n}\absb{\piff e_\three^1}
			+\tau^2\abs{e_\two ^1}
			+\tau^2\brat{\tau+\sigma_\four^n}\abs{\piff e_\two ^1}\Big)
			\quad\text{for $3\le n\leq N$.}
		\end{align*}
	\end{theorem}
	Theorem \ref{thm: error estimation of BDF2-DC3-DC4} together with the upper bounds in \eqref{ieq: bound decaying factors} 
	suggests that the BDF2-DC3-DC4 scheme is mesh-robustly fourth-order convergent on arbitrary time meshes if the starting solutions $v_\two ^1$, $v_\three^1$, $v_\four^1$ and $v_\four^2$ are accurate enough such that
	\begin{align*} 
		\abs{e_\four^{2}}
		+\tau\abs{\piff e_\four^2}+\tau\abs{e_\three^{1}}+\tau\brat{\tau+2^{-n}n}\absb{\piff e_\three^1}
		+\tau^2\abs{e_\two ^1}+\tau^2\brat{\tau+2^{-n}n^2}\abs{\piff e_\two ^1}=O(\tau^4).		
	\end{align*}
	Obviously, the error $\abs{e_\four^{2}}+\tau\abs{\piff e_\four^2}$ is logically determined by the first-level solver although the solution error of $v_\four^1$ is not explicitly involved in our error estimate. In practical applications, we always use the same starting solver to compute the two starting solutions $v_\four^1$ and $v_\four^2$.
	
	According to Remark \ref{remark: BDF2-DC3 starting}, the second-order L-stable RK method is suggested to start the BDF2 scheme and BDF2-DC3 scheme \eqref{eq: BDF2-DC3 scheme}. To achieve the fourth-order accuracy of the BDF2-DC3-DC4 scheme \eqref{eq: BDF2-DC3-DC4 scheme}, Theorem \ref{thm: error estimation of BDF2-DC3-DC4} says that the third-order RK method is required for the starting solutions $v_\four^1$ and $v_\four^2$.
	
	\begin{remark}\label{remark: BDF2-DC3-DC4 starting}	
	As done in Remark \ref{remark: BDF2-DC3 starting}, we list nine different cases for different starting groups (hereafter,  the triplet "\textbf{Scheme A+Scheme B+Scheme C}" represents that the starting solutions $v_\two ^1$, $v_\three^1$ and $v_\four^1 (v_\four^2)$ are generated by \textbf{Scheme A}, \textbf{Scheme B} and \textbf{Scheme C}, respectively):
		\begin{description}
			\item[(1)] Choose the RK3 method to start  the BDF2-DC3-DC4 scheme \eqref{eq: BDF2-DC3-DC4 scheme}
			\begin{description}
				\item[(RK2+RK2+RK3)]  $v_\two ^n=O(\tau^2)$, $v_\three^n=O(\tau^3)$ and $v_\four^n=O(\tau^4)$;
				\item[(BDF1+RK2+RK3)] $v_\two ^n=O(\tau^2)$, $v_\three^n=O(\tau^3)$ and  $v_\four^n=O(\tau^4)$;
				\item[(BDF1+BDF1+RK3)] $v_\two ^n=O(\tau^2)$,  $v_\three^n=O(\tau^2)$ and  $v_\four^n=O(\tau^3)$.
			\end{description}
			\item[(2)] Choose the RK2 method to start  the BDF2-DC3-DC4 scheme \eqref{eq: BDF2-DC3-DC4 scheme}
			\begin{description}
				\item[(RK2+RK2+RK2)]  
				$v_\two ^n=O(\tau^2)$,  $v_\three^n=O(\tau^3)$ and  $v_\four^n=O(\tau^3)$;
				\item[(BDF1+RK2+RK2)] $v_\two ^n=O(\tau^2)$, $v_\three^n=O(\tau^3)$ and  $v_\four^n=O(\tau^3)$;
				\item[(BDF1+BDF1+RK2)] $v_\two ^n=O(\tau^2)$, $v_\three^n=O(\tau^2)$ and  $v_\four^n=O(\tau^3)$.
			\end{description}
			\item[(3)] Choose the BDF1 method to start  the BDF2-DC3-DC4 scheme \eqref{eq: BDF2-DC3-DC4 scheme}
			\begin{description}
				\item[(RK2+RK2+BDF1)]  
				$v_\two ^n=O(\tau^2)$, $v_\three^n=O(\tau^3)$ and  $v_\four^n=O(\tau^2)$;
				\item[(BDF1+RK2+BDF1)] 	$v_\two ^n=O(\tau^2)$, $v_\three^n=O(\tau^3)$ and  $v_\four^n=O(\tau^2)$;
				\item[(BDF1+BDF1+BDF1)] $v_\two ^n=O(\tau^2)$, 	$v_\three^n=O(\tau^2)$ and  $v_\four^n=O(\tau^2)$;
			\end{description}
		\end{description}
		 It is observed that the own starting values $v_\four^1$ and $v_\four^2$ will have a significant impact on the accuracy of the BDF2-DC3-DC4 method; while the BDF2-DC3-DC4 method has no aftereffect, which means that the low-order starting methods that are accurate enough (can generate the desired accuracy) to start the BDF2 and BDF2-DC3 schemes will not cause a loss in the accuracy of the BDF2-DC3-DC4 method.  	
	\end{remark}

	\section{Extension to the BDF2-DC4 method}
	\setcounter{equation}{0}
	The BDF2-DC methods in the above two sections can improve one-order accuracy by one-step deferred correction process. One may wonder whether one-step deferred correction can improve two-order accuracy without sacrificing the mesh-robust stability. 
	
	In more detail, the BDF2 scheme is used to calculate the second-order solution $v_\two ^n$ and the fourth-order correction $\mathcal{D}_{2,4} f(t_n,v_\two ^n)$ is added in the BDF2 code to generate the fourth-order BDF2-DC4 solution $v_{\mathbf{4}'}^{n}$. We have the following numerical scheme
	\begin{align}
		\textbf{BDF2-DC4:}&\quad 
		D_2 v_{\mathbf{4}'}^{n} + \mathcal{D}_{2,4} f(t_n,v_\two ^n)   
		= f(t_n,v_{\mathbf{4}'}^n)\quad\text{for $n\geq 3$}, \label{eq: BDF2-DC4 scheme}	
	\end{align}
	where the correction operator $\mathcal{D}_{2,4}$ is defined by \eqref{de: D24}. To be more clear, we list the following figure to show how to startup the BDF2-DC4 scheme \eqref{eq: BDF2-DC4 scheme}
	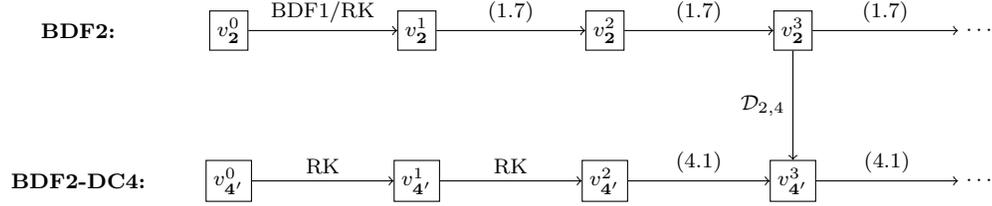
\begin{figure}[htb!]
		\centering
		\scriptsize
		\tikzstyle{format}=[rectangle,draw,thin,fill=white]
		\begin{tikzpicture}		
			\node[below of=BDF2,node distance=20mm] (BDF2){\textbf{BDF2:}};
			\node[format,right of=BDF2,node distance=20mm] (v330){$v_\two ^0$};
			\node[format,right of=v330,node distance=25mm] (v331){$v_\two ^1$};
			\node[format,right of=v331,node distance=25mm] (v332){$v_\two ^2$};
			\node[format,right of=v332,node distance=25mm] (v333){$v_\two ^3$};
			\node[right of=v333,node distance=25mm] (v334){$\cdots$};
			
			\node[below of=BDF2,node distance=20mm] (BDF2-DC4){\textbf{BDF2-DC4:}};
			\node[format,right of=BDF2-DC4,node distance=20mm] (v40){$v_{\mathbf{4}'}^0$};
			\node[format,right of=v40,node distance=25mm] (v41){$v_{\mathbf{4}'}^1$};
			\node[format,right of=v41,node distance=25mm] (v42){$v_{\mathbf{4}'}^2$};
			\node[format,right of=v42,node distance=25mm] (v43){$v_{\mathbf{4'}}^3$};
			\node[right of=v43,node distance=25mm] (v44){$\cdots$};	
			
			\draw[->](v330)--node[above]{BDF1/RK}(v331);
			\draw[->](v331)--node[above]{(1.7)}(v332);
			\draw[->](v332)--node[above]{(1.7)}(v333);
			\draw[->](v333)--node[above]{(1.7)}(v334);
			\draw[->](v333)--node[left]{$\mathcal{D}_{2,4}$}(v43);
			
			\draw[->] (v40)--node[above]{RK}(v41);
			\draw[->](v41)--node[above]{RK}(v42);
			\draw[->](v42)--node[above]{\eqref{eq: BDF2-DC4 scheme}}(v43);
			\draw[->](v43)--node[above]{\eqref{eq: BDF2-DC4 scheme}}(v44);
		\end{tikzpicture}
		\caption{The starting procedure (BDF1 or RK) of 
			BDF2-DC4 scheme \eqref{eq: BDF2-DC4 scheme}.}
		\label{figure: BDF2-DC4 initial process}
	\end{figure}
	Assume that $\bar{v}_{\mathbf{4}'}^{n}$ solves the following equation with a bounded sequence $\{\varepsilon_{\mathbf{4'}}^n\}$
	\begin{align*}
	\textbf{BDF2-DC4:}&\quad 
		D_2 \bar{v}_{\mathbf{4}'}^{n} + \mathcal{D}_{2,4} f(t_n,\bar{v}_\two ^n)= f(t_n,\bar{v}_{\mathbf{4}'}^n)+\varepsilon_{\mathbf{4}'}^n
		\quad \text{for $n\geq 3$}.
	\end{align*}
	The perturbed error
	$\tilde{v}_{\mathbf{4}'}^n:=\bar{v}_{\mathbf{4}'}^n-v_{\mathbf{4}'}^n$ satisfies the perturbed error system
	\begin{align}
		\textbf{BDF2-DC4:}&\quad 
		D_2 \tilde{v}_{\mathbf{4}'}^{n}= f(t_n,\bar{v}_{\mathbf{4}'}^n)-f(t_n,v_{\mathbf{4}'}^n)+\mathcal{D}_{2,4} f(t_n,v_\two ^n)-\mathcal{D}_{2,4} f(t_n,\bar{v}_\two ^n)
		+\varepsilon_{\mathbf{4}'}^n.\label{eq: BDF2-DC4 perturbed equation}
	\end{align}
	
	As done before, we consider the convolution form of \eqref{eq: BDF2-DC4 perturbed equation}.
	Multiplying both sides of the perturbed error equation \eqref{eq: BDF2-DC4 perturbed equation} by the DOC kernels
	$\vartheta_{n-i}^{(n)}$, and summing $i$ from 3 to $n$, we apply the identity \eqref{eq: BDF2-DOC identity} to get
	\begin{align}\label{eq: BDF2-DC4 DOC error equation}
		\piff \tilde{v}_{\mathbf{4}'}^n
		=-\mathcal{I}_{2}^n[\tilde{v}_{\mathbf{4}'}]+\sum\limits_{i=3}^n\vartheta_{n-i}^{(n)}
		\kbraB{f(t_i,\bar{v}_{\mathbf{4}'}^i)-f(t_i,v_{\mathbf{4}'}^i)
			+\mathcal{D}_{2,4}f(t_i,v_\two ^i)-\mathcal{D}_{2,4}f(t_i,\bar{v}_\two ^i)
			+\varepsilon_{\mathbf{4}'}^i}
	\end{align}
	for $3\leq n\leq N$, where $\mathcal{I}_{2}^n[\tilde{v}_{\mathbf{4}'}]$ represents the starting effect on the BDF2-DC4 solution error,
	\begin{align}\label{eq: initial effect BDF2-DC4}
		\mathcal{I}_{2}^n[\tilde{v}_{\mathbf{4}'}]
		:=&\,\sum_{i=3}^n\vartheta_{n-i}^{(n)}
		\sum_{j=1}^2d^{(i)}_{i-j}\piff \tilde{v}_{\mathbf{4}'}^j
		=\vartheta_{n-3}^{(n)}d^{(3)}_{1}\piff \tilde{v}_{\mathbf{4}'}^2
		\quad\text{for $3\le n\le N$.}
	\end{align}
	Following the derivation of \eqref{eq: D24 error integral form}, one can bound the perturbed error of  $\mathcal{D}_{2,4}\,f(t_i,v_\two ^i)$ by 
	\begin{align}\label{eq: D24 perturb error BDF2-DC4}
		\absb{\mathcal{D}_{2,4}\,f(t_i,v_\two ^i)-\mathcal{D}_{2,4}\,f(t_i,\bar{v}_\two ^i)}
		&\leq\absb{\mathcal{D}_{2,3}\,f(t_i,v_\two ^i)-\mathcal{D}_{2,3}\,f(t_i,\bar{v}_\two ^i)}
		+\frac{\tau_i}{6}\big|\hat{F}^{i}_\two +\hat{K}^{i}_\two 
		- \hat{F}^{i-1}_\two -\hat{K}^{i-1}_\two \big|\nonumber\\
		&\quad
		+\frac{\tau_ir_{i-1}(r_i+1)}{6(r_{i-1}+1)}\big|\hat{F}^{i-1}_\two +\hat{K}^{i-1}_\two 
		-\hat{F}^{i-2}_\two -\hat{K}^{i-2}_\two \big|,
	\end{align}
	where 
	\begin{align*}
		\hat{F}^{i}_\two :=&\,\int_0^1f'\big(\hat{q}_\two ^{i}\big)\partial_\tau\tilde{v}_\two ^i\zd s_1,\qquad\qquad\qquad\,\,\, \hat{q}_\two ^i:=(1-s_1)v_\two ^i+s_1\bar{v}_\two ^i,\\
		\hat{K}^{i}_\two :=&\,\int_0^1\int_0^1f''\big(\hat{g}_\two ^{i}\big)\piff \hat{q}_\two ^{i}\tilde{v}_\two ^{i-1}\zd s_1\zd s_2,\quad\hat{g}_\two ^i:=(1-s_2)\hat{q}_\two ^{i-1}+s_2\hat{q}_\two ^{i}.
	\end{align*} 
	By following the proof of Theorem \ref{thm: BDF2-DC3-DC4 stability}, one has the following stability result.
	\begin{theorem}\label{thm: BDF2-DC4 stability}
		Assume $\tau\le 1/(4L_f)$ and $f\in C^2$, the solution of \eqref{eq: BDF2-DC4 perturbed equation} satisfies
		\begin{align*}
			\abs{\tilde{v}_{\mathbf{4}'}^{n}}\leq4\exp(4L_ft_{n-1})
			\Big(\abs{\tilde{v}_{\mathbf{4}'}^{2}}
			+\tau\abs{\piff \tilde{v}_{\mathbf{4}'}^2}
			+t_n\max_{3\leq i\leq n}\Big\{\abs{\varepsilon_{\mathbf{4}'}^i}
			+\absb{\mathcal{D}_{2,4}f(t_i,v_\two ^i)-\mathcal{D}_{2,4}f(t_i,\bar{v}_\two ^i)}\Big\}\Big)
		\end{align*}
		for $3\le n\le N$. Thus the BDF2-DC4 scheme \eqref{eq: BDF2-DC4 scheme} is mesh-robustly stable.
	\end{theorem}
	
	To process the analysis, we consider the error effect of the correction term $\mathcal{D}_{2,4} f(t_n,v_\two ^n)$.
	\begin{lemma}\label{lem: truncation error of BDF2-DC4}
		Assume that $v(t)\in C^3((0,T])$ and $f\in C^2$. If $\tau\le 1/(4L_f)$, it holds that
		\begin{align*}
			\sum\limits_{i=3}^n\vartheta_{n-i}^{(n)}\absb{\mathcal{D}_{2,4} f(t_i,v_\two ^i)- \mathcal{D}_{2,4} f(t_i,v(t_i))}\leq Cr_{\max}\tau\bra{\abs{e_\two ^1}
				+\brat{\tau+\sigma_\three^n}\abs{\piff e_\two ^1}+t_n\tau^2},\;\;3\le n\leq N.
		\end{align*}	
	\end{lemma}
	\begin{proof}By following the proof of Lemma \ref{lem: truncation error of DC4}, it is not difficult to find that	
		\begin{align*}
			\abs{\mathcal{D}_{2,4}\,f(t_i,v_\two ^i)-\mathcal{D}_{2,4}\,f(t_i,v(t_i))}\le&\, \absb{\mathcal{D}_{2,3}\,f(t_i,v_\two ^i)-\mathcal{D}_{2,3}\,f(t_i,v(t_i))}\\
			&\,+Cr_{\max}\tau_i\sum_{\ell=i-2}^i\absb{\partial_\tau e_\two ^{\ell}}+Cr_{\max}\tau_i\sum_{\ell=i-3}^{i-1}\absb{e_\two ^{\ell}},\quad i\ge3.
		\end{align*}
		The proof of Lemma \ref{lem: truncation error of DC3} gives
		\begin{align*}
			&\absb{\mathcal{D}_{2,3}\,f(t_i,v_\two ^i)-\mathcal{D}_{2,3}\,f(t_i,v(t_i))}		
			\leq  C\tau_i\braB{\abs{\partial_\tau e_\two ^{i}}+\abs{\partial_\tau e_\two ^{i-1}}
				+\abs{e_\two ^{i-1}}+\abs{e_\two ^{i-2}}},\quad i\ge3.
		\end{align*}
		Thus it follows from Theorem \ref{thm: error estimation of BDF2}
		and Lemma \ref{lemma: BDF2 time derivative error} that 
		\begin{align*}
			&\absb{\mathcal{D}_{2,4}\,f(t_i,v_\two ^i)-\mathcal{D}_{2,4}\,f(t_i,v(t_i))}\\		
			&\leq Cr_{\max}\tau_i\sum_{\ell=i-2}^i\absb{\partial_\tau e_\two ^{\ell}}+Cr_{\max}\tau_i\sum_{\ell=i-3}^{i-1}\absb{e_\two ^{\ell}}\\
			&\leq Cr_{\max}\tau
			\bra{\abs{e_\two ^{1}}+(\tau+\sigma_\two ^i)\absb{\piff e_\two ^1}+t_i\tau^2}+Cr_{\max}\tau
			\braB{\abs{e_\two ^{1}}+\tau\absb{\piff e_\two ^1}+t_i\tau^2}\\
			&\leq Cr_{\max}\tau\bra{\abs{e_\two ^1}+\brat{\tau+\sigma_\two ^i}\abs{\piff e_\two ^1}+t_i\tau^2}.
		\end{align*}
		Recalling the definitions in  \eqref{def: decaying factors}, we apply Lemma \ref{lem: BDF2 orthogonal formula} to find that
		\begin{align*}
			\sum\limits_{i=3}^n\vartheta_{n-i}^{(n)}\absb{\mathcal{D}_{2,4} f(t_i,v_\two ^i)
				- \mathcal{D}_{2,4} f(t_i,v(t_i))}\leq Cr_{\max}\tau\bra{\abs{e_\two ^1}
				+\brat{\tau+\sigma_\three^n}\abs{\piff e_\two ^1}+t_n\tau^2}.
		\end{align*}	
		It completes the proof.	
	\end{proof}
	
	Now we present the convergence of BDF2-DC4 scheme \eqref{eq: BDF2-DC4 scheme}. Let numerical errors $e_{\mathbf{4}'}^0:=0$ and  $e_{\mathbf{4}'}^i:=v(t_i)-v_{\mathbf{4}'}^i$ for $1\le i\le N$.
	Subtracting \eqref{con: ODE DC equation} from \eqref{eq: BDF2-DC4 scheme}, one gets the following error equation of BDF2-DC4 scheme
	\begin{align}\label{eq: BDF2-DC4 error equation}
		D_2 e_{\mathbf{4}'}^{i} = f(t_i,v(t_i))-f(t_i,v_{\mathbf{4}'}^i)+\mathcal{D}_{2,4} f(t_i,v_\mathbf{2}^i)- \mathcal{D}_{2,4} f(t_i,v(t_i))+R_\four^i\quad\text{for $i\geq 3.$}
	\end{align}
	Multiplying both sides of \eqref{eq: BDF2-DC4 error equation} by the DOC kernels
	$\vartheta_{n-i}^{(n)}$, and summing $i$ from 3 to $n$, one gets
	\begin{align*}
		\piff e_{\mathbf{4}'}^n
		&=-\mathcal{I}_{2}^n[e_{\mathbf{4}'}]+\sum\limits_{i=3}^n\vartheta_{n-i}^{(n)}
		\braB{f(t_i,v(t_i))-f(t_i,v_{\mathbf{4}'}^i)+R_\four^i}\nonumber\\
		&\quad + \sum\limits_{i=3}^n\vartheta_{n-i}^{(n)}
		\braB{\mathcal{D}_{2,4}f(t_i,v_\two ^i)-\mathcal{D}_{2,4}f(t_i,v(t_i))}\quad \text{for $n\geq 3$},
	\end{align*}
	where the starting effect term $\mathcal{I}_{2}^n[e_{\mathbf{4}'}]$ is defined via \eqref{eq: initial effect BDF2-DC4}.
	Then, by using this error equation
	and Lemma \ref{lem: truncation error of BDF2-DC4}, 
	one can follow the proof of Theorem \ref{thm: BDF2-DC3-DC4 stability} to prove the following theorem.
	\begin{theorem}\label{thm: error estimation of BDF2-DC4}
		Let the solution $v\in C^5((0,T])$ and $v_\two ^n$ be the discrete solution  of the BDF2 scheme \eqref{eq: DC2 scheme}. If $\tau\le 1/(4L_f)$ and $f\in C^3$, the numerical solution $v_{\mathbf{4}’}^n$ of the BDF2-DC4 scheme \eqref{eq: BDF2-DC4 scheme} is convergent, 
		\begin{align*}
			\abs{e_{\mathbf{4}'}^{n}}\leq Cr_{\max}
			\Big(\abs{e_{\mathbf{4}'}^{2}}
			+\tau\abs{\piff e_{\mathbf{4}'}^2}+\tau\abs{e_\two ^1}
			+\tau\brat{\tau+\sigma_\three^n}\abs{\piff e_\two ^1}
			+t_n\tau^4\Big)
			\quad\text{for $3\le n\leq N$.}
		\end{align*}
	\end{theorem}
	Theorem \ref{thm: error estimation of BDF2-DC4} together with the upper bounds in \eqref{ieq: bound decaying factors} 
	suggests that the BDF2-DC4 scheme is mesh-robustly third-order convergent on arbitrary time meshes if the starting solutions $v_\two ^1$,  $v_{\mathbf{4}'}^1$ and $v_{\mathbf{4}'}^2$ are accurate enough such that
	\begin{align*} 
		\abs{e_{\mathbf{4}'}^{2}}
		+\tau\abs{\piff e_{\mathbf{4}'}^2}
		+\tau\abs{e_\two ^1}+\tau\brat{\tau+2^{-n}n}\abs{\piff e_\two ^1}=O(\tau^3).		
	\end{align*}
	
	Obviously, to achieve the fourth-order accuracy, the third-order RK method is required to compute the starting solutions $v_{\mathbf{4}'}^1$ and $v_{\mathbf{4}'}^2$. Moreover, the error estimate in Lemma \ref{lem: truncation error of BDF2-DC4} for the deferred correction $\mathcal{D}_{2,4}\,f(t_i,v_\two ^i)$  should be improved by imposing certain condition on the time meshes. If the step ratios $r_k$ vary slowly such that $\abs{r_{k+1}-r_{k}}\le C\tau$ and the starting solution $v_\two ^1$ is fourth-order accurate (via a third-order RK method) such that $\abs{\piff e_\two ^2-\piff e_\two ^1}\le C\tau^3$, the proof of \cite[Theorem 3.5]{BourgaultGaron:2022} yields
	\begin{align*}
		\absb{\mathcal{D}_{2,4} f(t_i,v_\two ^i)- \mathcal{D}_{2,4} f(t_i,v(t_i))}\leq C\tau^4,\quad3\le i\leq N.
	\end{align*}
	Then the result of Lemma \ref{lem: truncation error of BDF2-DC4} can be updated as follows,
	\begin{align*}
		\sum\limits_{i=3}^n\vartheta_{n-i}^{(n)}\absb{\mathcal{D}_{2,4} f(t_i,v_\two ^i)- \mathcal{D}_{2,4} f(t_i,v(t_i))}\leq C\tau^4,\quad3\le n\leq N.
	\end{align*}
	In such case, we have the following result.
	\begin{theorem}\label{thm: fourth-order error of BDF2-DC4}
		Let the solution $v\in C^5((0,T])$, $f\in C^3$ and $v_\two ^n$ be the discrete solution of the BDF2 scheme \eqref{eq: DC2 scheme} with a fourth-order starting value $v_\two ^1$. If $\tau\le 1/(4L_f)$ and $\abs{r_{k+1}-r_{k}}\le C\tau$ for $k\ge2$, the solution $v_{\mathbf{4}’}^n$ of the BDF2-DC4 scheme \eqref{eq: BDF2-DC4 scheme} is convergent in the sense that, 
		\begin{align*}
			\abs{e_{\mathbf{4}'}^{n}}\leq C
			\bra{\abs{e_{\mathbf{4}'}^{2}}
			+\tau\abs{\piff e_{\mathbf{4}'}^2}+t_n\tau^4}
			\quad\text{for $3\le n\leq N$.}
		\end{align*}
	\end{theorem}
	To end this section, we present a further comment on Theorems \ref{thm: error estimation of BDF2-DC4} and \ref{thm: fourth-order error of BDF2-DC4}. Although the BDF2-DC4 method has a loss of accuracy when the time-step sizes change rapidly, 
	this method would be a useful candidate in simulating some dissipative problems such as the gradient flows \cite{LiaoJiZhang:2020pfc,LiaoJiWangZhang:2021,LiaoKang2023,LiaoTangZhou:2020bdf2}.
	Actually, when the solution varies slowly (typically, in the long-time coarsening dynamics approaching the steady state), the uniform and quasi-uniform meshes with large time-step sizes would be preferred to resolve the numerical behavior. In simulating such slow dynamics, the BDF2-DC4 method is fourth-order accurate according 
	to Theorem \ref{thm: fourth-order error of BDF2-DC4} since the adjacent step ratios $r_k$ are
	always close to 1. In the fast-varying (high gradient) domains or in the transition regions between the
	slow-varying and fast-varying domains, the solution varies rapidly and a low-order but mesh-robust method with small time-step sizes would be preferred to capture the complex numerical behavior. In simulating such fast dynamics, the BDF2-DC4 method is mesh-robustly third-order accurate according 
	to Theorem \ref{thm: error estimation of BDF2-DC4} because the  step ratios $r_k$ have some relatively large fluctuations.

	\section{Numerical experiments}
	\setcounter{equation}{0}
	This section presents extensive numerical examples, including tests on different meshes, effects of different startings, behavior for stiff problem and adaptive time-stepping, to examine the stability and convergence of the BDF2-DC methods 
	\eqref{eq: DC2 scheme}-\eqref{eq: BDF2-DC3-DC4 scheme} and \eqref{eq: BDF2-DC4 scheme}. 
	We record the  error $e(N):=\max_{1\leq{n}\leq{N}}|v(t_n)-v^n|$ and compute the numerical order of convergence by
	$\text{Order}\approx\log\big(e(N)/e(2N)\big)/\log\big(\tau(N)/\tau(2N)\big),$
	where $\tau(N)$ represents the maximum step size for total $N$ subintervals. 
	\subsection{Tests on different meshes}
	\begin{example}\cite{BourgaultGaron:2022}\label{example1}
		Consider the model $v_t=v\cos t$ for $0<t<T$ such that 
		$v(t)=e^{\sin t}$. 
	\end{example}
	\begin{table}[htb!]
		\begin{center}
			\caption{The adjacent time-step ratios $r_{k}$ under two meshes for Example \ref{example1} with $T=10\pi$. 
				\label{table: rk for example1}}
			\vspace*{0.3pt}
			\def\temptablewidth{1\textwidth}
			{\rule{\temptablewidth}{0.5pt}}
			\begin{tabular*}{\temptablewidth}{@{\extracolsep{\fill}}ccccccc}
				&\multirow{2}{*}{$N$}  &\multicolumn{2}{c}{G-mesh}   &\multicolumn{2}{c}{R-mesh}  \\
				\cline{3-4}          \cline{5-6}          
				&  &$\tau/\tau_{1}(\gamma=2)$ &$\tau/\tau_{1}(\gamma=3)$  &$r_{max}$ &$N_1$  \\
				\midrule         
				&5120 &1.02E+04 &7.86E+07 &4.48E+03 &1035  \\
				&10240 &2.05E+04 &3.15E+08 &6.66E+03 &2116  \\	
				&20480 &4.10E+04 &1.26E+09 &3.66E+04 &4222  \\
			\end{tabular*}
			{\rule{\temptablewidth}{0.5pt}}
		\end{center}
	\end{table}	
	
	This subsection examines the BDF2-DC methods for Example \ref{example1} on different meshes: 
	\begin{itemize}[itemindent=-5mm]
		\item\textbf{G-mesh}: $t_k=T(k/N)^\gamma$ with some grading parameters $\gamma>1$ such that $r_{max}=2^\gamma-1$. To show the defect of Definition \ref{Def:stability1} (zero-stability), see the disscusions in Remark \ref{rem: zero-stable condition}, we record the stability prefactor $\tau/\tau_1=\tau_N/\tau_1=N^{\gamma}-(N-1)^{\gamma}$ in  Table \ref{table: rk for example1} for two different parameters. 
		\item \textbf{R-mesh}: random steps $\tau_k=\epsilon_kT/\sum_{\ell=1}^N\epsilon_\ell$ with uniformly distributed random numbers $\epsilon_\ell\in (0,1)$. Table \ref{table: rk for example1} lists the maximum step ratio $r_{\max}$ and and the number (denote $N_1$) of time levels with the step ratio $r_k> 1+\sqrt{2}$ on three R-meshes for different $N$. As seen, $r_{\max}$ is very large and there are many of (about $20\%$ in our tests) 
		step-ratios which are greater than the classical stability limit. 
		\item \textbf{S-mesh}: $t_k=T\cdot3^{k-N}$ with a fixed step-ratio $r_k=3$ for $2\le k\le N$.
	\end{itemize}	
	
	\begin{table}[htb!]
		\begin{center}
			\caption{Numerical results of BDF2-DC methods for Example \ref{example1} on G-meshes with $T=10\pi$. 
				\label{table: graded meshes for example1}} \vspace*{0.3pt}
			\def\temptablewidth{1.0\textwidth}
			{\rule{\temptablewidth}{0.5pt}}
			\begin{tabular*}{\temptablewidth}{@{\extracolsep{\fill}}ccccccccc}
				\multirow{2}{*}{Scheme} &\multirow{2}{*}{$N$}  &\multicolumn{2}{c}{$\gamma=2$} &\multirow{2}{*}{CPU time} &\multicolumn{2}{c}{$\gamma=3$} &\multirow{2}{*}{CPU time} \\
				\cline{3-4}          \cline{6-7}         
				&  &$e(N)$ &Order&  &$e(N)$ &Order& \\
				\midrule  
				&	5120	&	3.79E-05	&	-	&	1.94E-02	&	8.46E-05	&	-	&	2.14E-02\\
				BDF2  &	10240	&	9.45E-06	&	2.01	&	2.87E-02	&	2.11E-05	&	2.00 	&	2.72E-02\\
				&	20480	&	2.36E-06	&	2.00 	&	6.62E-02	&	5.26E-06	&	2.00 &	4.41E-02\\
				\hline                 
				&	5120	&	9.18E-08	&	-	&	3.55E-02	&	1.82E-07	&	-	&	3.40E-02\\
				BDF2-DC3 &	10240	&	1.15E-08	&	3.00	&	4.44E-02	&	2.28E-08	&	3.00	&	4.42E-02\\
				&	20480	&	1.44E-09	&	3.00	&	9.46E-02	&	2.87E-09	&	2.99	&	8.14E-02\\
				\hline    
				&	5120	&	2.62E-07	&	-	&	2.36E-02	&	4.79E-07	&	-	&	1.95E-02\\
				BDF3  	&	10240	&	3.37E-08	&	2.96	&	3.34E-02	&	6.44E-08	&	2.90	&	2.86E-02\\
				&	20480	&	4.27E-09	&	2.98	&	5.22E-02	&	8.36E-09	&	2.95	&	5.05E-02\\
				\hline                                               
				&	5120	&	2.15E-09	&	-	&	5.55E-02	&	1.05E-08	&	-	&	5.07E-02\\
				BDF2-DC3-DC4 &	10240	&	1.46E-10	&	3.88	&	7.56E-02	&	7.38E-10	&	3.83	&	7.03E-02\\
				&	20480	&	9.46E-12	&	3.95	&	1.34E-01	&	4.87E-11	&	3.92	&	1.27E-01\\
				\hline   
				&	5120	&	8.83E-09	&	-	&	2.38E-02	&	4.26E-08	&	-	&	2.16E-02\\
				BDF4		&	10240	&	6.11E-10	&	3.85	&	3.18E-02	&	3.07E-09	&	3.80	&	3.31E-02\\
				&	20480	&	3.99E-11	&	3.94	&	5.64E-02	&	2.04E-10	&	3.91	&	5.61E-02\\
			\end{tabular*}
			{\rule{\temptablewidth}{0.5pt}}
		\end{center}
	\end{table}	
	
	\begin{table}[htb!]
		\begin{center}
			\caption{Numerical results of BDF2-DC methods for Example \ref{example1} on R-meshes with $T=10\pi$. 
				\label{table: random meshes for example1}}
			\vspace*{0.3pt}
			\def\temptablewidth{1\textwidth}
			{\rule{\temptablewidth}{0.5pt}}
			\begin{tabular*}{\temptablewidth}{@{\extracolsep{\fill}}cccccccc}
				&\multirow{2}{*}{$N$}  &\multicolumn{2}{c}{BDF2}&\multicolumn{2}{c}{BDF2-DC3}   &\multicolumn{2}{c}{BDF2-DC3-DC4}  \\
				\cline{3-4}          \cline{5-6}          \cline{7-8}   
				&     &$e(N)$     &Order  &$e(N)$  &Order &$e(N)$ &Order     \\
				\midrule         
				&5120  &5.26E-06 &-  &1.23E-07 &- &2.88E-08 & \\
				&10240 &1.21E-06 &2.11 &1.09E-08 &3.50 &1.63E-09 &4.14 \\	
				&20480 &1.85E-07 &2.73 &1.36E-09 &2.99 &7.84E-11 &4.36 \\
			\end{tabular*}
			{\rule{\temptablewidth}{0.5pt}}
		\end{center}
	\end{table}	
	
		To remove the impact of startup errors, the solutions in  
		Tables \ref{table: graded meshes for example1}-\ref{table: BDF2-DC4 scheme for example1} are always generated by using the exact starting values $v_\two ^1=v(t_1), v_\three^1=v(t_1), v_\four^1=v(t_1)$ and $v_\four^2=v(t_2)$. 
		Table \ref{table: graded meshes for example1}
		lists the numerical results and CPU time of five different methods on G-meshes with two grading parameters $\gamma=2$ and $\gamma=3$.
		It can be observed that the numerical solutions remain stable and convergent with full accuracy although the ratio $\tau/\tau_1$ is quite large, cf. the second and third columns in Table \ref{table: rk for example1}. In this sense, also see Remark \ref{rem: zero-stable condition}, Definition 2 and the stability results in Lemma \ref{lem: BDF2 stability}, Theorems \ref{thm: BDF2-DC3 stability} and \ref{thm: BDF2-DC3-DC4 stability} are practically reasonable. From the fifth and eighth columns in Table \ref{table: graded meshes for example1}, we see that, as expected, the BDF2-DC methods really cost more CPU times than the standard BDF3 and BDF4 schemes. This is practically reasonable since the DC approaches always need the solutions of two
		and three implicit systems per step to achieve the same accuracy. In fact, these additional computational cost are the price we pay for theoretically improving the rigid step-ratio conditions (cf. \cite[Section III.5]{HairerNrsettWanner1993} or the recent development in \cite{LiLiao:2022}) for the stability and convergence of variable-step BDF3 and BDF4 schemes.
	
		Table \ref{table: random meshes for example1} records the solution errors of the BDF2, BDF2-DC3 and BDF2-DC3-DC4 method on R-meshes.  
		The results suggest that the variable-step BDF2-DC methods are mesh-robustly stable and convergent. They support the stability estimates in Lemma \ref{lem: BDF2 stability}, Theorems \ref{thm: BDF2-DC3 stability} and \ref{thm: BDF2-DC3-DC4 stability}. The experimental orders of convergence also meet our predictions in
		Theorems \ref{thm: error estimation of BDF2}, \ref{thm: error estimation of BDF2-DC3} and \ref{thm: error estimation of BDF2-DC3-DC4}.
	
	\begin{table}[htb!]
		\begin{center}
			\caption{Numerical results of BDF2-DC methods for Example \ref{example1} on S-mesh with $T=1$. 
				\label{table: special meshes for example1}}
			\vspace*{0.3pt}
			\def\temptablewidth{1\textwidth}
			{\rule{\temptablewidth}{0.5pt}}
			\begin{tabular*}{\temptablewidth}{@{\extracolsep{\fill}}cccccccccc}
				&\multirow{2}{*}{$N$} &\multirow{2}{*}{$r_{max}$}&\multirow{2}{*}{$N_1$} &\multicolumn{2}{c}{BDF2}&\multicolumn{2}{c}{BDF2-DC3}   &\multicolumn{2}{c}{BDF2-DC3-DC4}  \\
				\cline{5-6}          \cline{7-8}          \cline{9-10}   
				&  &&   &$e(N)$     &Order  &$e(N)$  &Order &$e(N)$ &Order     \\
				\midrule         
				&10 &3 &8  &1.40E-01 &- &2.05E-02 &- &2.02E-03 &- \\
				&20 &3 &18 &1.40E-01 &0.00 &2.05E-02 &0.00 &2.02E-03 &0.00 \\	
				&40 &3 &38 &1.40E-01 &0.00 &2.05E-02 &0.00 &2.02E-03 &0.00 \\
			\end{tabular*}
			{\rule{\temptablewidth}{0.5pt}}
		\end{center}
	\end{table}	
	
	\begin{table}[htb!]
		\begin{center}
			\caption{Numerical results of BDF2-DC4 scheme for Example \ref{example1} with $T=1$. 
				\label{table: BDF2-DC4 scheme for example1}}
			\vspace*{0.3pt}
			\def\temptablewidth{1\textwidth}
			{\rule{\temptablewidth}{0.5pt}}
			\begin{tabular*}{\temptablewidth}{@{\extracolsep{\fill}}cccccccc}
				&\multirow{2}{*}{$N$}  &\multicolumn{2}{c}{G-meshes ($\gamma=2$)}   &\multicolumn{2}{c}{R-meshes}  &\multicolumn{2}{c}{S-mesh}\\
				\cline{3-4}          \cline{5-6}      \cline{7-8}    
				&   &$e(N)$  &Order   &$e(N)$ &Order  &$e(N)$ &Order  \\
				\midrule         
				&10 &2.59E-04 &- &2.00E-04 &- &1.53E-02 &- \\
				&20 &2.22E-05 &3.69 &2.79E-05 &2.74 &1.53E-02 &0.00 \\	
				&40 &1.60E-06 &3.86 &4.19E-06 &2.84 &1.53E-02 &0.00 \\
			\end{tabular*}
			{\rule{\temptablewidth}{0.5pt}}
		\end{center}
	\end{table}	
	
	As a special case in examining the zero-stability, we run the BDF2, BDF2-DC3 and BDF2-DC3-DC4 methods on the S-mesh, which has $r_k=3$ for all $k\ge2$. As seen in Table \ref{table: special meshes for example1}, the numerical solutions remain bounded (stable) although they are always not convergent as $N$ increases (because the maximum step-size $\tau_N=2/3$ is unchanged for different $N$). Again, they say that the classical zero-stability condition  $0<r_k< 1+\sqrt{2}$ is not necessary.
		
	 We also run the BDF2-DC4 scheme \eqref{eq: BDF2-DC4 scheme} on the above three time meshes and list the BDF2-DC4 solution errors in Table \ref{table: BDF2-DC4 scheme for example1}. The numerical results support the stability estimate in Theorem \ref{thm: BDF2-DC4 stability}. It is observed that, as predicted by Theorems \ref{thm: error estimation of BDF2-DC4} and \ref{thm: fourth-order error of BDF2-DC4}, the BDF2-DC4 scheme \eqref{eq: BDF2-DC4 scheme} is fourth-order accurate on the G-meshes, third-order accurate on the R-meshes, and zero-order accurate on the S-mesh. Actually, the graded meshes asymptotically  satisfy the mesh condition $\abs{r_{k+1}-r_{k}}\le C\tau$  of Theorem \ref{thm: fourth-order error of BDF2-DC4} as  the level index $k$ is properly large.
	 
	\subsection{Effects of different starting solutions}
	\begin{table}[htb!]
		\begin{center} 
			\caption{The accuracy of BDF2-DC3-DC4 scheme initiated by RK3 for Example \ref{example1} with $T=10\pi$.
				\label{table: RK3 uniform meshes for example1}} \vspace*{0.3pt}
			\def\temptablewidth{1.0\textwidth}
			{\rule{\temptablewidth}{0.5pt}}
			\begin{tabular*}{\temptablewidth}{@{\extracolsep{\fill}}ccccccccc}
				\multirow{2}{*}{Starting schemes} &\multirow{2}{*}{$N$} &\multicolumn{2}{c}{BDF2 } &\multicolumn{2}{c}{BDF2-DC3} &\multicolumn{2}{c}{BDF2-DC3-DC4}  \\
				\cline{3-4}          \cline{5-6}          \cline{7-8}
				& &$e(N)$ &Order &$e(N)$ &Order  &$e(N)$ &Order \\
				\midrule      
				&	1280	&	6.27E-04	&		&	1.11E-03	&		&	8.10E-07	&	-\\
				\textbf{BDF1+BDF1+RK3}&	2560	&	1.54E-04	&	2.02	&	2.92E-04	&	1.93	&	1.20E-07	&	2.76\\
				&	5120	&	3.82E-05	&	2.01	&	7.49E-05	&	1.96	&	1.62E-08	&	2.89\\
				\hline                          
				&	1280	&	1.38E-03	&		&	1.07E-04	&		&	6.59E-07	&	-\\
				\textbf{BDF1+RK2+RK3}&	2560	&	3.48E-04	&	1.98	&	1.34E-05	&	3.00	&	4.11E-08	&	4.00\\
				&	5120	&	8.76E-05	&	1.99	&	1.67E-06	&	3.00	&	2.57E-09	&	4.00\\
				\hline 
				&	1280	&	9.86E-04	&		&	1.13E-04	&		&	7.36E-07	&	-\\
				\textbf{RK2+RK2+RK3}&	2560	&	2.48E-04	&	1.99	&	1.42E-05	&	3.00	&	4.59E-08	&	4.00\\
				&	5120	&	6.23E-05	&	1.99	&	1.78E-06	&	3.00	&	2.87E-09	&	4.00\\
			\end{tabular*}
			{\rule{\temptablewidth}{0.5pt}}
		\end{center}
	\end{table}	
	\begin{table}[htb!]
		\begin{center} 
			\caption{ The accuracy of BDF2-DC3-DC4 scheme initiated by RK2 for Example \ref{example1} with $T=10\pi$.
				\label{table: RK2 uniform meshes for example1}} \vspace*{0.3pt}
			\def\temptablewidth{1.0\textwidth}
			{\rule{\temptablewidth}{0.5pt}}
			\begin{tabular*}{\temptablewidth}{@{\extracolsep{\fill}}ccccccccc}
				\multirow{2}{*}{Starting schemes} &\multirow{2}{*}{$N$} &\multicolumn{2}{c}{BDF2 } &\multicolumn{2}{c}{BDF2-DC3} &\multicolumn{2}{c}{BDF2-DC3-DC4}  \\
				\cline{3-4}          \cline{5-6}          \cline{7-8}
				& &$e(N)$ &Order &$e(N)$ &Order  &$e(N)$ &Order \\
				\midrule      
				&	1280	&	6.27E-04	&		&	1.11E-03	&		&	3.85E-06	&	-\\
				\textbf{BDF1+BDF1+RK2}&	2560	&	1.54E-04	&	2.02	&	2.92E-04	&	1.93	&	5.00E-07	&	2.95\\
				&	5120	&	3.82E-05	&	2.01	&	7.49E-05	&	1.96	&	6.38E-08	&	2.97\\
				\hline                            
				&	1280	&	1.38E-03	&		&	1.07E-04	&		&	2.50E-06	&	-\\
				\textbf{BDF1+RK2+RK2}&	2560	&	3.48E-04	&	1.98	&	1.34E-05	&	3.00	&	3.45E-07	&	2.86\\
				&	5120	&	8.76E-05	&	1.99	&	1.67E-06	&	3.00	&	4.53E-08	&	2.93\\
				\hline 
				&	1280	&	9.86E-04	&		&	1.13E-04	&		&	2.41E-06	&	-\\
				\textbf{RK2+RK2+RK2}&	2560	&	2.48E-04	&	1.99	&	1.42E-05	&	3.00	&	3.40E-07	&	2.83\\
				&	5120	&	6.23E-05	&	1.99	&	1.78E-06	&	3.00	&	4.50E-08	&	2.92\\
			\end{tabular*}
			{\rule{\temptablewidth}{0.5pt}}
		\end{center}
	\end{table}	
	
	\begin{table}[htb!]
		\begin{center} 
			\caption{The accuracy of BDF2-DC3-DC4 scheme initiated by BDF1 for Example \ref{example1} with $T=10\pi$.
				\label{table: BDF1 uniform meshes for example1}} \vspace*{0.3pt}
			\def\temptablewidth{1.0\textwidth}
			{\rule{\temptablewidth}{0.5pt}}
			\begin{tabular*}{\temptablewidth}{@{\extracolsep{\fill}}ccccccccc}
				\multirow{2}{*}{Starting schemes} &\multirow{2}{*}{$N$} &\multicolumn{2}{c}{BDF2 } &\multicolumn{2}{c}{BDF2-DC3} &\multicolumn{2}{c}{BDF2-DC3-DC4}  \\
				\cline{3-4}          \cline{5-6}          \cline{7-8}
				& &$e(N)$ &Order &$e(N)$ &Order  &$e(N)$ &Order \\
				\midrule           
				&	1280	&	6.27E-04	&	-	&	1.11E-03	&	-	&	2.00E-03	&	-\\
				\textbf{BDF1+BDF1+BDF1}&	2560	&	1.54E-04	&	2.02	&	2.92E-04	&	1.93	&	5.05E-04	&	1.98\\
				&	5120	&	3.82E-05	&	2.01	&	7.49E-05	&	1.96	&	1.27E-04	&	1.99\\
				\hline                            
				&	1280	&	1.38E-03	&	-	&	1.07E-04	&	-	&	1.99E-03	&	-\\
				\textbf{BDF1+RK2+BDF1}&	2560	&	3.48E-04	&	1.98	&	1.34E-05	&	3.00	&	5.05E-04	&	1.98\\
				&	5120	&	8.76E-05	&	1.99	&	1.67E-06	&	3.00	&	1.27E-04	&	1.99\\
				\hline  
				&	1280	&	9.86E-04	&		&	1.13E-04	&	-	&	1.99E-03	&	-\\
				\textbf{RK2+RK2+BDF1}&	2560	&	2.48E-04	&	1.99	&	1.42E-05	&	3.00	&	5.05E-04	&	1.98\\
				&	5120	&	6.23E-05	&	1.99	&	1.78E-06	&	3.00	&	1.27E-04	&	1.99\\
			\end{tabular*}
			{\rule{\temptablewidth}{0.5pt}}
		\end{center}
	\end{table}	
	
	To investigate the effects  of the starting approximations, we will compute the starting solutions $v_\two ^1$, $v_\three^1$ and $v_\four^1 (v_\four^2)$ by the two-stage q-order SDIRK (in short, RKq) methods for $q=2$ or $q=3$ as follows,
		\begin{align*}
			\text{RK2}:\quad\begin{array}{c|cc}
				\frac{2-\sqrt{2}}{2} &\frac{2-\sqrt{2}}{2}  &   \\ 
				1 &1-\frac{2-\sqrt{2}}{2}  &\frac{2-\sqrt{2}}{2} \\
				\hline 
				&1-\frac{2-\sqrt{2}}{2}  &\frac{2-\sqrt{2}}{2}  \\ 
			\end{array}\qquad\quad		
			\text{RK3}:\quad	\begin{array}{c|cc}
				\frac{3+\sqrt{3}}{6} &\frac{3+\sqrt{3}}{6}  &   \\ 
				1-\frac{3+\sqrt{3}}{6} &1-\frac{3+\sqrt{3}}{3}  &\frac{3+\sqrt{3}}{6} \\
				\hline 
				&\frac12   &\frac12 \\ 
			\end{array}
	\end{align*}
	Without special declarations, we always use the triplet ``\textbf{Scheme A+Scheme B+Scheme C}", see Remarks \ref{remark: BDF2-DC3 starting} and \ref{remark: BDF2-DC3-DC4 starting}, to represent that  $v_\two ^1$, $v_\three^1$ and $v_\four^1 (v_\four^2)$ are generated by \textbf{Scheme A}, \textbf{Scheme B} and \textbf{Scheme C}, respectively.  We examine the impact of startup errors by considering nine different cases for three different starting groups, listed in Tables \ref{table: RK3 uniform meshes for example1}-\ref{table: BDF1 uniform meshes for example1}, corresponding to the RK3, RK2 and BDF1 starting schemes, respectively, for the BDF2-DC3-DC4 scheme \eqref{eq: BDF2-DC3-DC4 scheme} on the uniform mesh. 
	
	As seen from the last two columns in Tables \ref{table: RK3 uniform meshes for example1}-\ref{table: BDF1 uniform meshes for example1}, when the BDF2 and BDF2-DC3 methods achieve their full accuracy, the BDF2-DC3-DC4 scheme \eqref{eq: BDF2-DC3-DC4 scheme} can generate the desired fourth-order accuracy only if a third-order scheme is applied to start the procedure; while a loss of accuracy is obvious when the RK2 and BDF1 starting schemes are used.
	
	Similarly, from the fifth and sixth columns in Tables \ref{table: RK3 uniform meshes for example1}-\ref{table: BDF1 uniform meshes for example1}, when the BDF2 method achieves the second-order accuracy, the BDF2-DC3 scheme \eqref{eq: BDF2-DC3 scheme} can generate the desired third-order accurate solution only if a second-order scheme is applied to start the procedure; while a loss of accuracy is obvious when the BDF1 starting scheme is used. These observations support the theoretical predictions in Remarks \ref{remark: BDF2-DC3 starting} and \ref{remark: BDF2-DC3-DC4 starting}. They suggest that our error estimates in Theorems \ref{thm: error estimation of BDF2}, \ref{thm: error estimation of BDF2-DC3} and \ref{thm: error estimation of BDF2-DC3-DC4} seem sharp especially on the numerical effects  of the starting approximations.
	
	\subsection{Numerical behaviors for stiff problem}
	\begin{example}\label{example2}\cite{BourgaultGaron:2022} Consider the following stiff differential system
		\begin{align*}
			& \frac{d \mathbf{u}}{d t}=\left(\begin{array}{ccc}
					-1 & 1 & 100 \\
					0 & 0 & 100 \\
					0 & -100 & 0
				\end{array}\right)\,\mathbf{u}\quad \text{with}\quad\mathbf{u}(0)=\left(\begin{array}{c}
				2 \\
				1 \\
				1
			\end{array}\right)
		\end{align*}
		on the interval $0\leq t\leq5$. The exact solution is given by the following equation
		\begin{align*}
			\mathbf{u}(t)=e^{-t}\left(\begin{array}{l}
				1 \\
				0 \\
				0
			\end{array}\right)+\cos (100 t)\left(\begin{array}{l}
				1 \\
				1 \\
				1
			\end{array}\right)+\sin (100 t)\left(\begin{array}{c}
				1 \\
				1 \\
				-1
			\end{array}\right).
		\end{align*}
	\end{example}
	
	\begin{table}[htb!]
		\begin{center}
			\caption{
				Numerical results of BDF2-DC methods for Example \ref{example2} on G-meshes.\label{table: graded meshes for example2}} \vspace*{0.3pt}
			\def\temptablewidth{1.0\textwidth}
			{\rule{\temptablewidth}{0.5pt}}
			\begin{tabular*}{\temptablewidth}{@{\extracolsep{\fill}}ccccccc}
				\multirow{2}{*}{Scheme} &\multirow{2}{*}{$N$}  &\multicolumn{2}{c}{$\gamma=2$}  &\multicolumn{2}{c}{$\gamma=3$}  \\
				\cline{3-4}          \cline{5-6}         
				&  &$e(N)$ &Order &$e(N)$ &Order \\
				\midrule                       
				&	100000	&	1.17E-02	&		&	2.26E-02	&	-\\
				&	200000	&	2.93E-03	&	2.00	&	5.65E-03	&	2.00\\
				BDF2&	400000	&	7.33E-04	&	2.00	&	1.41E-03	&	2.00\\
				&	800000	&	1.83E-04	&	2.00	&	3.53E-04	&	2.00\\
				&	1600000	&	4.58E-05	&	2.00	&	8.83E-05	&	2.00\\
				\hline
				&	100000	&	7.12E-05	&		&	2.43E-04	&	-\\
				&	200000	&	5.90E-06	&	3.59	&	1.93E-05	&	3.66\\
				BDF2-DC3&	400000	&	5.51E-07	&	3.42	&	1.71E-06	&	3.49\\
				&	800000	&	5.71E-08	&	3.27	&	1.71E-07	&	3.33\\
				&	1600000	&	6.41E-09	&	3.16	&	1.86E-08	&	3.20\\
				\hline
				&	100000	&	3.31E-07	&		&	1.88E-06	&	-\\
				&	200000	&	1.17E-08	&	4.82	&	5.79E-08	&	5.02\\
				BDF2-DC3-DC4&	400000	&	5.49E-10	&	4.41	&	2.47E-09	&	4.55\\
				&	800000	&	3.03E-11	&	4.18	&	1.28E-10	&	4.26\\
				&	1600000	&	7.76E-12	&	1.97	&	6.18E-12	&	4.38\\
			\end{tabular*}
			{\rule{\temptablewidth}{0.5pt}}
		\end{center}
	\end{table}	
	
	\begin{table}[htb!]
		\begin{center} 
			\caption{The accuracy of BDF2-DC3-DC4 scheme initiated by different methods for Example \ref{example2}.
				\label{table: uniform meshes for example2}} \vspace*{0.3pt}
			\def\temptablewidth{1.0\textwidth}
			{\rule{\temptablewidth}{0.5pt}}
			\begin{tabular*}{\temptablewidth}{@{\extracolsep{\fill}}ccccccccc}
				\multirow{2}{*}{Starting schemes} &\multirow{2}{*}{$N$} &\multicolumn{2}{c}{BDF2 } &\multicolumn{2}{c}{BDF2-DC3} &\multicolumn{2}{c}{BDF2-DC3-DC4}  \\
				\cline{3-4}          \cline{5-6}          \cline{7-8}
				& &$e(N)$ &Order &$e(N)$ &Order  &$e(N)$ &Order \\
				\midrule       
				&	1000	&	2.04E+00	&		&	3.26E+00	&		&	5.88E+00	&	-\\
				\textbf{BDF1+RK2+BDF1}&	4000	&	2.43E+00	&	-0.13 	&	3.32E+00	&	-0.01 	&	2.97E+00	&	0.49 \\
				&	16000	&	2.29E-01	&	1.70 	&	2.02E-02	&	3.68 	&	2.02E-03	&	5.26 \\
				&	64000	&	1.43E-02	&	2.00 	&	1.00E-04	&	3.83 	&	1.08E-04	&	2.11 \\
				\hline
				&	1000	&	2.04E+00	&		&	3.26E+00	&		&	5.85E+00	&	-\\
				\textbf{BDF1+RK2+RK2}&	4000	&	2.43E+00	&	-0.13 	&	3.32E+00	&	-0.01 	&	2.96E+00	&	0.49 \\
				&	16000	&	2.29E-01	&	1.70 	&	2.02E-02	&	3.68 	&	1.17E-03	&	5.65 \\
				&	64000	&	1.43E-02	&	2.00 	&	1.00E-04	&	3.83 	&	4.44E-07	&	5.68 \\
				\hline
				&	1000	&	2.04E+00	&		&	3.26E+00	&		&	5.85E+00	&	-\\
				\textbf{BDF1+RK2+RK3}&	4000	&	2.43E+00	&	-0.13 	&	3.32E+00	&	-0.01 	&	2.96E+00	&	0.49 \\
				&	16000	&	2.29E-01	&	1.70 	&	2.02E-02	&	3.68 	&	1.18E-03	&	5.65 \\ 
				&	64000	&	1.43E-02	&	2.00 	&	1.00E-04	&	3.83 	&	5.12E-07	&	5.58 \\
			\end{tabular*}
			{\rule{\temptablewidth}{0.5pt}}
		\end{center}
	\end{table}	
	
	Since BDF type methods are intended for stiff problems, we examine the numerical behaviors by running the BDF2 and BDF2-DC methods for Example \ref{example2}.
	Table \ref{table: graded meshes for example2}  records the solution errors on the G-meshes with the exact starting values  $v_\two ^1=v(t_1), v_\three^1=v(t_1), v_\four^1=v(t_1)$ and $v_\four^2=v(t_2)$.  The numerical results suggest that the variable-step BDF2, BDF2-DC3 and BDF2-DC3-DC4 methods are mesh-robustly stable and convergent.
	They support the stability estimates in Lemma \ref{lem: BDF2 stability}, Theorems \ref{thm: BDF2-DC3 stability} and \ref{thm: BDF2-DC3-DC4 stability}, and partially support the error estimates in
	Theorems \ref{thm: error estimation of BDF2}, \ref{thm: error estimation of BDF2-DC3} and \ref{thm: error estimation of BDF2-DC3-DC4}. 
	
	The numerical effects of the starting approximations for the BDF2, BDF2-DC3 and BDF2-DC3-DC4 methods with the uniform time-step are presentd in Table \ref{table: uniform meshes for example2}, which includes three starting triplet groups: \textbf{BDF1+RK2+RK3},   \textbf{BDF1+RK2+RK2} and  \textbf{BDF1+RK2+BDF1}. The numerical behaviors are similar to those in \cite[Table 10]{BourgaultGaron:2022}. As expected, a very small time-step size would be necessary when we use the uniform mesh to capture the multiscale behaviors in the stiff system.
	Experimentally, they support the sharp error estimates in Theorems \ref{thm: error estimation of BDF2}, \ref{thm: error estimation of BDF2-DC3} and \ref{thm: error estimation of BDF2-DC3-DC4}.  
	Actually, when the BDF2 method achieves the second-order accuracy, the BDF2-DC3 scheme \eqref{eq: BDF2-DC3 scheme} with a second-order starting scheme is  third-order convergent; and when the BDF2 and BDF2-DC3 methods achieve their full accuracy, the BDF2-DC3-DC4 scheme \eqref{eq: BDF2-DC3-DC4 scheme} with a third-order starting scheme is fourth-order convergent.
	
	\subsection{Adaptive time-stepping}
	To resolve the multiscale behaviors in stiff problems, one would need certain adaptive time-stepping approach, that is, to choose a small step size when the solution changes rapidly, and adopt a large step size when the solution changes slowly. We now examine the adaptive time-stepping strategy of the BDF-DC3 method for the following Allen-Cahn-type nonlinear model.
	\begin{example}\label{example3}\cite{StuartHumphries:1998,XuXu2023}
		Consider the model $v_t=v-v^3$, which has a monotone solution
		$$v(t)=\frac{v_0}{\sqrt{e^{-2t}+v_0^2(1-e^{-2t})}}\quad\rightarrow\quad \frac{v_0}{\abs{v_0}}\qquad\text{as $t\rightarrow+\infty$}.$$
		A simple fixed-point iteration algorithm with the termination error $10^{-12}$ 
		is applied to solve the resulting nonlinear equation at each time level. We compute the solution $v(t)$ by an adaptive time-stepping strategy.
	\end{example}
	
	Notice that, in the step-by-step solution procedure from the BDF2 method \eqref{eq: DC2 scheme} to the BDF2-DC3 method \eqref{eq: BDF2-DC3 scheme}, the latter can increase the accuracy of BDF2 solution from second to third order if the solution is smooth while the BDF2-DC3 solution may retain the second-order accuracy if the solution varies rapidly. It presents an immediate error estimator, such as the relative error $e^n_{23}:=|v_\three^{n}-v_\two ^{n}|/|v_\two ^{n}|$, for certain adaptive time-stepping strategy to choose the optimal time-step size $\tau_n$. Similarly, in the solution procedures from  \eqref{eq: BDF2-DC3 scheme} to  \eqref{eq: BDF2-DC3-DC4 scheme} and from \eqref{eq: DC2 scheme} to \eqref{eq: BDF2-DC4 scheme}, one can use the immediate error estimators  $e^n_{34}:=|v_\four^{n}-v_\three^{n}|/|v_\three^{n}|$ and $e^n_{24'}:=|v_{\mathbf{4'}}^{n}-v_\two ^{n}|/|v_\two ^{n}|$, respectively, to choose the optimal time-step size. 
	
	\begin{table*}[htb!]
		\centering
		\caption{CPU time (in seconds) and total time levels of the BDF2-DC3 scheme.}\label{table: CPU for example3}
		\begin{tabular}{c|cc|cc}
			\hline
			\multirow{2}{*}{$T$}& \multicolumn{2}{c|}{adaptive strategy}& \multicolumn{2}{c}{uniform mesh} \\
			\cline{2-5}
			&CPU time &Time levels &CPU time &Time levels \\
			\hline
			100    &0.23 &$10^3$      &1.77  &$10^5$	\\ 
			\hline
			1000  &0.40	&$10^4$	&13.68  &$10^6$	\\ 
			\hline
			10000  &3.92	&$10^5$    &142.78   &$10^7$ \\ 
			\hline
		\end{tabular}
	\end{table*}
	
	\begin{algorithm}[htb!]
		{\small 
		\SetAlgoLined
		\KwIn{the solution $v_\two ^{n}$ and step-size $\tau_{n-1}$.}
		\KwOut{the solution $v_\three^{n}$ and step-size $\tau_{n}$.}
		
		Calculate $v_\three^{n}$ by BDF2-DC3 scheme.
		
		Calculate relative error $e^n_{23}=|v_\three^{n}-v_\two ^{n}|/|v_\two ^{n}|$.
		
		\eIf{$e^n_{23}>tol$ }
		{Reject $v_\three^{n}$.\\		
			Recalculate step-size $\tau_{n}\leftarrow\max\{\tau_{\min},\tau_{ada}\}$; Goto 1.
		}
		{Accept $v_\three^{n}$.\\
			Update step-size $\tau_{n}\leftarrow\min\{\max\{\tau_{\min},\tau_{ada}\},\tau_{\max}\}$.\\
		}
	
		\caption{Adaptive time-stepping strategy for BDF2-DC3 method}
		\label{algorithm:adaptive BDF2-DC3}}
	\end{algorithm}
	
	\begin{figure*}
		\centering
		\subfigure[$v_0=-1.5$]{
			\begin{minipage}[b]{0.23\linewidth}
				\includegraphics[width=1\linewidth]{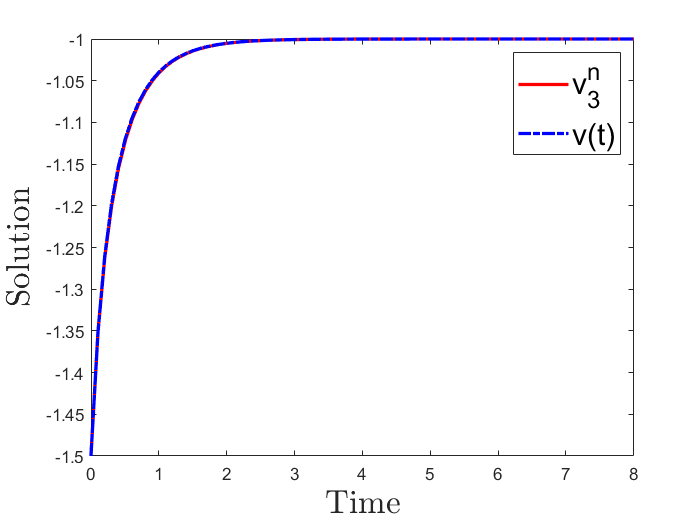}\vspace{4pt}
				\includegraphics[width=1\linewidth]{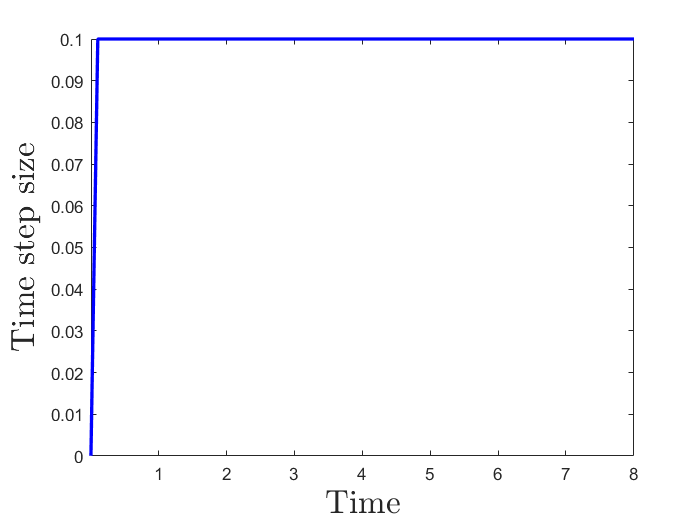}\vspace{4pt}
		\end{minipage}}
		\subfigure[$v_0=-0.5$]{
			\begin{minipage}[b]{0.23\linewidth}
				\includegraphics[width=1\linewidth]{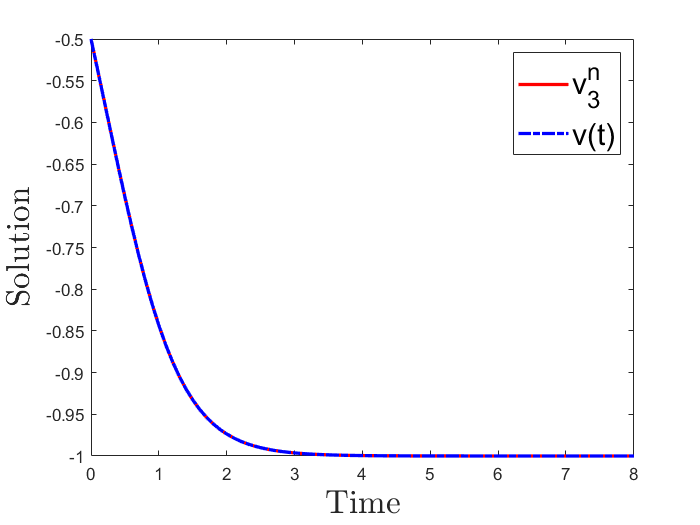}\vspace{4pt}
				\includegraphics[width=1\linewidth]{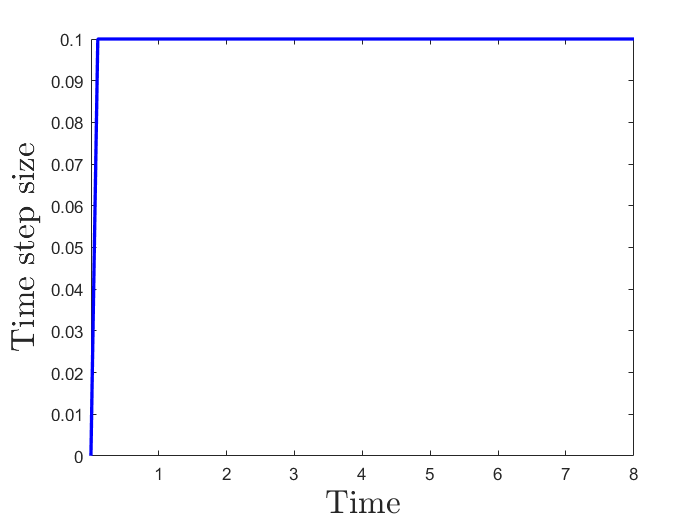}\vspace{4pt}
		\end{minipage}}
		\subfigure[$v_0=0.5$]{
			\begin{minipage}[b]{0.23\linewidth}
				\includegraphics[width=1\linewidth]{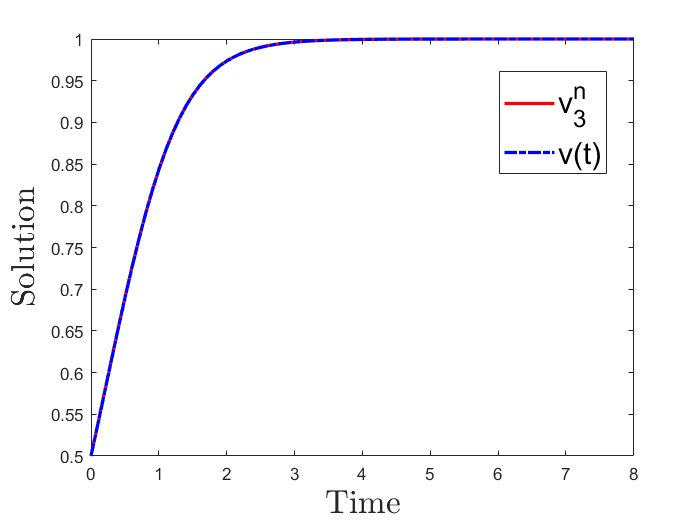}\vspace{4pt}
				\includegraphics[width=1\linewidth]{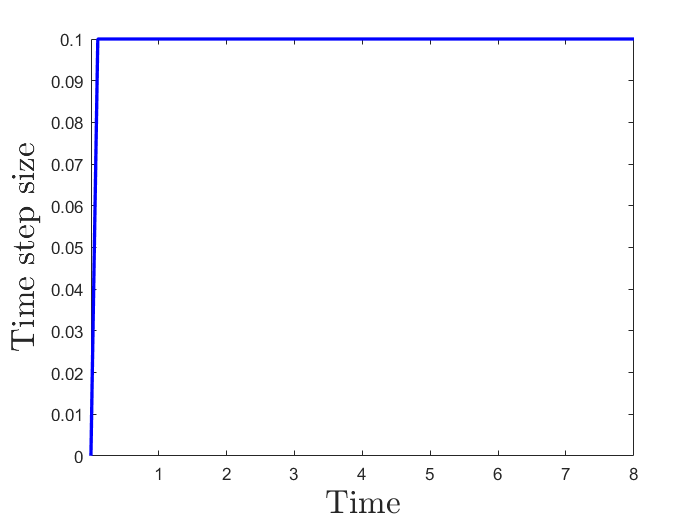}\vspace{4pt}
		\end{minipage}}
		\subfigure[$v_0=1.5$]{
			\begin{minipage}[b]{0.23\linewidth}
				\includegraphics[width=1\linewidth]{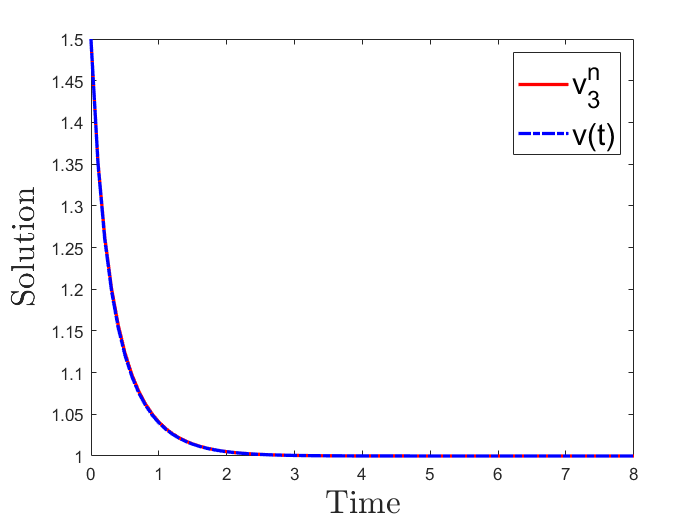}\vspace{4pt}
				\includegraphics[width=1\linewidth]{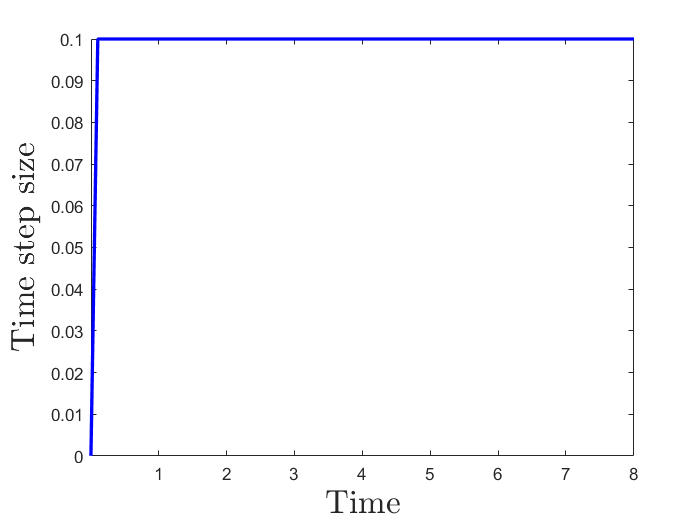}\vspace{4pt}
		\end{minipage}}
		\caption{Solution curves (up) and adaptive time-steps (down) for four initial values.}
		\label{figure:u0}
	\end{figure*}
	
	Here we adopt the adaptive time-stepping strategy from \cite{GomezHughes:2011,LiaoJiZhang:2020pfc} by using the relative error $e^n_{23}$ to determine the adaptive time-step size $\tau_{ada}$,
	$\tau_{ada}:=S_a\tau_n\sqrt{tol/e^n_{23}},$
	where $S_a$ is the safety parameter and $tol$ is the tolerance error. In Algorithm \ref{algorithm:adaptive BDF2-DC3} for the solution procedure from the BDF2 method \eqref{eq: DC2 scheme} to the BDF2-DC3 method \eqref{eq: BDF2-DC3 scheme}, we set $S_a=10^3$, $tol=10^{-1}$, the maximum step size $\tau_{max}=0.1$, the minimum step size $\tau_{min}=10^{-3}$ and the first level $t_1=\tau_{\min}$. 
	
	Table \ref{table: CPU for example3} compares the CPU times and total time levels in simulating Example \ref{example3} by using the uniform mesh and Algorithm \ref{algorithm:adaptive BDF2-DC3}. Since the solution of this example varies fast only in the early time period, one can use large time-steps in the long-time dynamics such that the adaptive time-stepping strategy is more computationally efficient (this advantage would be very valuable in solving the corresponding partial differential problems because the large-scale system of linear or nonlinear algebraic equations requires more CPU time at each time level). Also, Figure \ref{figure:u0} depicts the solution curves and the associated time-steps generated by Algorithm \ref{algorithm:adaptive BDF2-DC3} for four initial data $v_0=-1.5$, $-0.5$, $0.5$ and $1.5$. It seems that, at least for this example, the third-order BDF2-DC3 scheme maintains the monotonicity of continuous solution and converges to the correct steady state solution; however, the intrinsic mechanism is open to us and needs further study, cf. the recent work \cite{XuXu2023}.
	
	\section{Concluding remarks}
	By using the discrete orthogonal convolution kernels, we show that high-order BDF2-DC methods are mesh-robustly stable and convergent on arbitrary time grids and would be superior to the high-order BDF-k schemes \cite{LiLiao:2022,LiaoTangZhou2023,LiaoTangZhou2023BDF3parabolic} in the sense of the step-ratio condition. The present error estimates are novel and sharp in the sense that they clearly address the numerical effects of the starting approximations. We find that the BDF2-DC methods have no aftereffect
	in the case that one-step deferred correction process only improve one-order of accuracy. Applications to nonlinear parabolic equations are under current investigations and will be presented in separate reports.

\end{document}